\documentclass[11pt]{amsart}

\usepackage{amsthm}
\usepackage{breqn}
\usepackage{amsmath}
\usepackage{amsfonts}
\usepackage{amssymb}
\usepackage{tikz}
\usepackage{hyperref}
\usepackage{microtype} 
\usepackage{mathtools}
\usepackage{graphicx}

\newcommand{\R}{\mathbb{R}} 
\newcommand{\e}{\varepsilon}
\newcommand{\N}{\mathbb N}
\newcommand{\E}{\varphi}
\newcommand{\Z}{\mathcal{Z}_{n-1}}
\newcommand{\f}{{\bf f}}
\newcommand{\M}{{\bf M}}
\newcommand\Is{{\mathfrak{Is}}}

\newtheorem{theorem}{Theorem}[section]
\newtheorem*{theorem*}{Theorem}

\newtheorem{remarq}{Remark}

\title[\resizebox{4.6 in}{!}{Min-max for phase transitions and the existence of embedded minimal hypersurfaces}]{Min-max for phase transitions and the existence of embedded minimal hypersurfaces}
\author{Marco A. M. Guaraco}
\address{Instituto de Matem\'atica Pura e Aplicada (IMPA) \\ Estrada Dona Castorina 110 \\ 22460-320 Rio de Janeiro \\ Brazil}
\email{marcomen@impa.br}

\thanks{The author was partly supported by CAPES-Brazil and NSF Grant DMS-1104592.}

\begin{document}
 
 
\begin{abstract}

Strong parallels can be drawn between the theory of minimal hypersurfaces and the theory of phase transitions. Borrowing ideas from the former we extend recent results on the regularity of stable phase transition interfaces \cite{TonegawaWickramasekera} to the finite Morse index case. As an application we present a PDE-based proof of the celebrated theorem of Almgren-Pitts, on the existence of embedded minimal hypersurfaces in compact manifolds. We compare our results with other min-max theories; \cite{Pitts,SimonSmith}. \end{abstract}

\maketitle


\setcounter{tocdepth}{1}

\section{Introduction}

In \cite{AllenCahn}, Allen-Cahn introduced the semilinear parabolic equation \begin{align}\label{flow}\partial_t u -\Delta u + \frac{W'(u)}{\e^2}=0\end{align} as a mathematical model for the evolution of phase transition phenomena, the function $W$ being a double-well potential with unique global minima at $\pm1$. A typical example of such a non-linearity is $W(u)=\frac{1}{4}(1-u^2)^2$.

From a geometrical point of view solutions to equation (\ref{flow}) have a remarkable feature: roughly speaking, as $\e\to 0$ the level sets of $u$ concentrate around a hypersurface (called the \textit{limit-interface}) that is evolving in time under the action of the mean curvature flow (see for example \cite{Ilmanen,MizunoTonegawa,PisanteAllen13}). This suggests that for the particular case of stationary states of (\ref{flow}), i.e. solutions to the semilinear elliptic equation \begin{align}\label{eq:theequation}-\e\Delta u + \frac{W'(u)}{\e}=0,\end{align} the limit-interface should be a stationary point of the mean curvature flow, i.e. a minimal hypersurface.

Beginning with the works of L. Modica \cite{Modica} and P. Sternberg \cite{Sternberg}, this idea has been made precise in a variety of situations, and in the last decades the bridge between the theory of phase transitions and minimal hypersurfaces has been exploited extensively in both directions. We refer the reader to the surveys \cite{Pacard,Savin,TonegawaSurvey} and the references therein.

From a variational perspective solutions to equation (\ref{eq:theequation}) in a bounded open set $U\subset \R^n$, are critical points of the energy functional $$E_\e(u)=\int_U \e \frac{|\nabla u|^2}{2} + \frac{W(u)}{\e},$$ and variational properties such as stability or finite Morse index of a solution $u$ are defined as usual, i.e. with respect to the bilinear form corresponding to $E_\e''(u)(\cdot,\cdot)$, the second derivative of the energy in $H^1(U)$. 

Of special interest to us are the combined works of Hutchinson-Tonegawa \cite{HutchinsonTonegawa}, Tonegawa \cite{Tonegawa} and Tonegawa-Wickramasekera \cite{TonegawaWickramasekera}. Roughly speaking, they showed the following (see \ref{regularity})

\begin{theorem*}[\cite{HutchinsonTonegawa,Tonegawa,TonegawaWickramasekera}]
Let $U$ be a bounded open set in $\R^n$ and $u_{k}$ a sequence of solutions to \text{(\ref{eq:theequation})} in $U$, with $\e=\e_k\to0$. Assume that $\sup_U |u_k|$ and $E_{\e_k}(u_k)$ are bounded sequences. Then, as $\e_k\to 0$, the level sets of $u_k$ accumulate around a stationary integral varifold. If in addition the solutions are stable, then the limit-interface is a stable minimal hypersurface, smooth and embedded outside a set of Hausdorff dimension at most $n-8$.
\end{theorem*}

We extend the regularity statement to the case of solutions with finite Morse index on compact manifolds. The idea that regularity of general critical points can be obtained from the stable case, goes back to the work of Pitts \cite{Pitts} on the min-max construction of minimal hypersurfaces. We adapt these ideas (see the technical remark at the end of this section) to the phase transition context  to obtain (see \ref{morseindex})

\subsection*{Theorem A}{\it
Let $M$ be a $n$-dimensional closed Riemannian manifold and $u_{k}$ a sequence of solutions to \text{(\ref{eq:theequation})} in $M$, with $\e=\e_k\to0$. Assume that their Morse indices, $\sup_M |u_k|$ and $E_{\e_k}(u_k)$ are bounded sequences. Then, as $\e_k\to 0$, its level sets accumulate around a minimal hypersurface, smooth and embedded outside a set of Hausdorff dimension at most $n-8$.}

\

There is a natural family of solutions to which one can try to apply this result. In fact, since the constant functions $\pm1$ are the only global minimizers of the energy  $$E_\e(u)=\int_M \e \frac{|\nabla u|^2}{2} + \frac{W(u)}{\e},$$ we can expect to obtain other critical points by min-max methods. After checking a Palais-Smale condition, we use an extension of the mountain-pass theorem (see \cite{Ghoussoub}, Chapter 10) to show the existence of solutions $u_\e$ with Morse index at most 1. After dividing by twice the energy constant $\sigma=\int_{-1}^{1}\sqrt{W(s)/2}ds$, the energies $c_\e = E_\e(u_\e)$ of these solutions will converge, as $\e\to0$, to the area of the limit-interface given by Theorem A. More precisely, we prove (see \ref{phasetransitionminimal})

\subsection*{Theorem B}\label{phasetransitionminimal} \textit{In every $n$-dimensional closed Riemannian manifold there exists an integral varifold $V$ such that
\begin{enumerate}
\item [(i)] $\|V\|=\frac{1}{2\sigma}\liminf c_\e$;
\item [(ii)] $V$ is stationary in $M$;
\item [(iii)] $\mathcal{H}^{n-8+\gamma}(\operatorname{sing}(V))=0,$ for every $\gamma>0$;
\item [(iv)] $\operatorname{reg}(V)$ is an embedded minimal hypersurface.
\end{enumerate}
}

\ 

As a Corollary we obtain the celebrated 
\begin{theorem*}[Almgren, Pitts, Schoen-Simon]
Every $n$-dimensional closed Riemannian manifold contains a minimal hypersurface, smooth and embedded outside a set of Hausdorff dimension at most $n-8$.
\end{theorem*}
 
An important step in our construction is controlling the behavior of the energies $c_\e$ as $\e\to 0$. In order to guarantee the existence of a limit-interface these energies cannot explode or vanish. In this respect, we prove the following upper bound (see Section \ref{upperbound.section})
 
\subsection*{Proposition C}{ \it Let $\{\Sigma_t\}_{t\in[0,1]}$ be a sweepout given by isotopic deformations of the level sets of a Morse function. Then $$\frac{1}{2\sigma}\limsup c_\e \leq \max_{t\in[0,1]} \mathcal{H}^{n-1}(\Sigma_t).$$

} 

\

A similar result can be proven for sweepouts coming from Heegaard splittings, with essentially the same proof.

We also obtain lower bounds for $\liminf c_\e$ in two different ways. First, a standard application of the Sobolev space analogue of the isoperimetric inequality (that goes back to De Giorgi) gives a short proof of the fact that energies of the solution do not vanish as $\e \to 0$ (see Section \ref{bounds.section}).  This is important because it guarantees that the limit-interface is not trivial. Second, a sharper lower bound is obtained by constructing discrete sweepouts of currents, with arbitrarily small fineness and mass controlled by the energies of mountain-pass solutions to (\ref{eq:theequation}). In what follows $\Pi$ represents the fundamental homotopy class of one-parameter sweepouts containing those given by level sets of Morse functions. The relevant definitions can be found in Section \ref{almgrenpitts.section}.

\subsection*{Proposition D} {\it There is a non-trivial homotopy class $\Pi$ (see Propositon \ref{nontrivialclass}) such that $$0<{\bf L}(\Pi)\leq \frac{1}{2\sigma}\liminf c_\e.$$ In particular, if $V$ is the limit-interface of Theorem B, and $V_{AP}$ is the stationary varifold obtained after applying Almgren-Pitt's min-max to the class $\Pi$ (Theorem \ref{almgrenpitts}), then $$\|V_{AP} \| \leq \|V\|.$$}

When $n=3$, the upper bound from Proposition C implies that the area of the limit-interface is no bigger than the area of the surface obtained by Simon-Smith's continuous refinement of the min-max methods (presented by Colding-De Lellis in \cite{ColdingDelellis}). More precisely, we have  

\subsection*{Corollary E}{\it If $n=3$ and $V_{cont}$ is the minimal surface (with multiplicities) obtained with the continuous min-max methods of \cite{ColdingDelellis}, when applied to the saturated family of sweepouts generated by isotopic deformations of the level sets of a Morse function (or coming from a Heegaard splitting), then $$ \|V_{AP} \| \leq  \|V \| \leq \|V_{cont}\|.$$} 

In some special cases (e.g. if $Ric_M>0$, or in the absence of stable hypersurfaces) it has been proved that a minimal hypersurface of least area has index 1 and area $\|V_{AP}\|= \|V_{cont}\|$ (see \cite{MarquesNevesDuke,Rosenberg,Zhou,Zhou2015}). The formula above gives us the same index bound, when $n = 3$, for the minimal surface obtained by phase-transition methods.

Theorem B is an example of a one-parameter min-max construction. We note that our methods can be applied to families with more parameters whenever they satisfy the requirements of a mountain-pass lemma (see \cite{Ghoussoub}, Chapter 10). In recent years min-max constructions have been proved to be a powerful tool in the understanding of minimal hypersurfaces in different context. We refer the reader to  \cite{DelellisPellandini,Ketover,MarquesNevesDuke,Rosenberg,Zhou,Zhou2015} for applications to compact manifolds; \cite{RosenbergHyp,Montezuma} for applications to non-compact manifolds; \cite{Li} for applications to a free-boundary situation; and to Marques-Neves \cite{MarquesNevesWillmore,MarquesNevesInfinite} for the solution of the long-standing Willmore conjecture and the existence of infinite embedded minimal hypersurfaces in compact manifolds with positive Ricci curvature, respectively.

Finally, we briefly mention some results emphasizing the parallel between the theory of phase transitions and the earlier developed theory of minimal hypersurfaces. Regularity results for minimal hypersurfaces were first obtained for area minimizing currents in the solution of Plateau's problem by De Giorgi and Federer-Fleming. A regularity theory for stable minimal hypersurfaces was developed by Schoen-Simon-Yau \cite{SchoenSimonYau} (and later by Schoen-Simon \cite{SchoenSimon}), and used by Pitts to show the regularity of unstable minimal hypersurfaces in the context of Almgren's min-max theory of varifolds. In the same way, first regularity results for the convergence of phase transitions were obtained by De Giorgi's school at Pisa in the 70's, for energy minimizing solutions. Motivated by the work of Pitts, Padilla-Tonegawa \cite{TonegawaPadilla}, Hutchinson-Tonegawa \cite{HutchinsonTonegawa}, Tonegawa \cite{Tonegawa} and  Tonegawa-Wickramasekera \cite{TonegawaWickramasekera}, carried out a program to obtain a weak convergence theory for general (unstable) phase transitions and the regularity for the stable case. The present work is a natural continuation of these results.

\subsection*{Technical remark} It is worth mentioning some technicals points concerning the proofs of Theorems A and B in comparison with Pitts' technique. 

The notion of {\it almost minimizing} varifolds has a fundamental role in the work of Pitts. It allows one to use the estimates of Schoen-Simon-Yau \cite{SchoenSimonYau} (and later Schoen-Simon \cite{SchoenSimon}) for the construction of stable {\it replacements}, an essential step in obtaining the regularity of the limit varifold. In general, a min-max limit object is expected to have Morse index at most the number of parameters used for the construction.  We believe that one of the motivations behind Pitt's results is the observation that finite index implies stability in small annular regions of the ambient space (see Remark \ref{indexstable}). Although this is a feature of variational nature that holds with great generality, something like finiteness of the Morse index is very hard to check in his context (even for particular cases e.g. see \cite{MarquesNevesWillmore}). The introduction of the almost minimizing property by Pitts seems like an attempt to overcome this difficulty by proving the stability of the limit varifold in annular regions directly from the variational construction, without referring to any notion of index. 

In the phase-transitions context the situation is slightly different. We use a min-max lemma only in the proof of Theorem B, to guarantee the existence of solutions of the PDE (\ref{eq:theequation}). The notion of Morse index for these solutions is clear, and there are several results in the literature concerning index bounds (see \cite{Ghoussoub}, Chapter 10). In particular, the stability of the solutions in small annular regions follows from this, without the introduction of any almost minimizing property. On the other hand, the regularity of the stable limit-interfaces used in the proof of Theorem A, proven by Tonegawa-Wickramasekera in \cite{TonegawaWickramasekera}, depends on a non-trivial extension of the estimates of \cite{SchoenSimon} due to N. Wickramasekera \cite{Wickramasekera}.

\subsection*{Organization}

The content of this work is organized as follows.

In Section \ref{preliminaries}, we present some notation and review preliminaries from the theory of varifolds and the theory of stable minimal hypersurfaces.

In Section \ref{convergence}, we show that the regularity results from \cite{TonegawaWickramasekera}, for stable phase transition interfaces, can be extended to the finite Morse index case. 

In Section \ref{minmax.section}, using standard mountain-pass techniques for elliptic operators, we show the existence of solutions to equation ($\ref{eq:theequation}$) with index at most 1.

In Section \ref{limit-interface.section}, we combine the results from the previous sections to show the existence of embedded minimal hypersurfaces in compact manifolds. We also sketch the ideas that motivate the computations of the energy bounds in Sections \ref{bounds.section} and  \ref{upperbound.section}.

In Section \ref{bounds.section}, we prove a lower bound for the energy of mountain-pass solutions to equation (\ref{eq:theequation}), using a Sobolev space version of the isoperimetric inequality.

In Section \ref{upperbound.section}, we give upper bounds for the energy of mountain-pass solutions to equation (\ref{eq:theequation}). We show that the limit-interface has area no bigger than the width of level sets of Morse functions. Some technical details concerning distance functions that are used in this section are developed in one appendix (Section \ref{appendix.section}).

In Section \ref{almgrenpitts.section}, we show that the limit-interface has area no smaller than the minimal hypersurfaces obtained with Almgren-Pitts min-max theory. For $n=3$, we show that it has area no bigger than the surface obtained by the continuous version of Almgren-Pitts theory from \cite{ColdingDelellis}.

Finally, in another appendix (Section \ref{manifolds.section}), we comment on the extension to general manifolds of the results from \cite{TonegawaPadilla,HutchinsonTonegawa,Tonegawa,TonegawaWickramasekera}.

\subsection*{Acknowledgements}

This work is partially based on my Ph.D. thesis at IMPA. I am grateful to my advisor, Fernando Cod\'a Marques, for his constant encouragement and support. I am also grateful to the Mathematics Department of Princeton University for its hospitality. The first drafts of this work were written there while visiting during Fall of 2014.

\

\

\section{Notation and preliminaries}\label{preliminaries}

\subsection{Notation} Along this work we will use the following notation

\medskip

\begin{tabular}{lll}
$W$ && a double-well potential (see \ref{assumptions})\\
$\sigma$ && the energy constant $\sigma=\int_{-1}^{1}\sqrt{W(s)/2}ds$\\  
$\operatorname{Inj}(M)$ && the injectivity radius of $M$. \\
$An(x,\tau,\tau')$ && the annulus centered at $x$ with radii $0<\tau<\tau'$.\\
$\mathcal{AN}(x,r)$ && the set $\{ An(x,\tau,\tau') : 0<\tau<\tau'<r \}$. \\
$\mathcal{B}_r(x)$ && the ball centered at 0 with radius $r$ in $T_x M$.\\
$B(x,r)$ && the geodesic ball centered at $x$ with radius $r$.\\
$\Is(M)$ && the set of isotopies of $M$\\
$W^{1,p}(M)$ && is the Sobolev space of $L^p(M)$ functions \\
&& with weak derivatives also in $L^p(M)$.\\
$H^1(M)$ && the Sobolev space $W^{1,2}(M)$ \\
$d_K(\cdot)$ && the distance function from a closed set $K\subset M.$\\
\end{tabular}

\subsection{The theory of varifolds} 
In order to study general variational problems, Almgren introduced the notion of \textit{varifold} as a generalization of the concept of submanifolds. The theory of currents of Federer and Fleming, available at the time, was suitable for treating minimization problems, but for generalizing the min-max technique of Birkhoff a new notion was necessary. The reason for this is in part that in the theory of currents the area functional is only lower-semicontinous and some cancelation of mass may occur when dealing with limits that are saddle points.

Let $U$ be an open subset of a Riemannian manifold $M^n$. We denote by $G(U)$ the $(n-1)$-dimensional Grassmanian bundle of unoriented hyperplanes over $U$.

An \textit{$(n-1)$-varifold} in $U$, or simply a $\textit{varifold}$ for the extent of this work, is any nonnegative, finite Radon measure on $G(U)$. The space of varifolds is endowed with the topology of weak* convergence, so a sequence of varifolds $V_k$ converge to a varifold $V$ if for every $\varphi\in C_c(G(U))$ we have $$\int \varphi(x,\pi) dV_k(x,\pi) \longrightarrow  \int \varphi(x,\pi) dV_k(x,\pi).$$ 

We can associate a positive measure on $U$ to any varifold V, we call it the \textit{mass} of $V$, and is defined by the formula $$\int_U \varphi(x) d\|V\|(x) = \int_{G(U)} \varphi(x) dV(x,\pi).$$ We also refer to $\|V\|(U)$ as the \emph{mass of $V$ in $U$}. This measure generalizes the area functional.

Given any $(n-1)$-rectifiable set $\Sigma\subset U$ and a $\mathcal{H}^{n-1}$-measurable function $\theta:\Sigma \to \mathbb{Z}$, called the \textit{multiplicity}, we can associate to them a varifold $V_{\theta \Sigma}$ by $$V_{\theta \Sigma}(\varphi)=\int \varphi(x,T_x \Sigma) \theta(x) d\mathcal{H}^{n-1}(x),$$ for any $\varphi \in C_c(G(U))$. Any varifold obtained in that way is called an \textit{integer varifold}. When $\theta \equiv 1$ we simply write $V_\Sigma$.

By using the change of variables formula we can \textit{pushforward} a varifold $V$ in the presence of a diffeormophism $\psi: U\to U'$, by defining the varifold $$ \psi_* (V)(\varphi)=\int_U \phi(\psi(x),d_x\psi(\pi))|d_x\psi| dV(x,\pi),$$ where $|d_x\psi|$ is the Jacobian of $\psi$ at $x$.

Finally, we can use this pushforward to define notions of first and second variation for varifolds with respect to vector fields in $U$. Given a smooth vector field $X$ supported in $U$, denote $\psi(t)$ the associated flow (i.e. $d\psi(t) / dt =X$). We define the \textit{first} and \textit{second variations of $V$ with respect to $X$}, respectively, by $$[\delta V](X)= \frac{d}{dt}\|\psi(t)_* V \|(U)\bigg|_{t=0}$$ and $$[\delta^2 V](X)= \frac{d^2}{dt^2}\|\psi(t)_* V \|(U)\bigg|_{t=0}.$$ As in the theory of smooth manifolds, we say that a varifold is \emph{stationary} if $[\partial V](X)=0$, for any vector field $X$, and we call it \emph{stable}, if in addition $[\partial^2 V](X)\geq0$.

\subsection{Stable hypersurfaces}
Given an orientable open set $U\subset M$ and $\Gamma \subset M$ a closed subset of codimension 1, we say that $\Gamma$ satisfies (SS) in $U$, if  $\Gamma \cap U$ is a smooth embedded hypersurface outside a closed set $S$, with $\mathcal{H}^{n-3}(S) = 0$.

 The following theorem, taken from \cite{DelellisTasnady}, is a consequence of Schoen-Simon curvature estimates. 

\begin{theorem}\label{stablecompacity} Let $U$ be an orientable open subset of a manifold and $\{g^k\}$ and $\{\Gamma^k\}$, respectively, sequences of smooth metrics on $U$ and of hypersurfaces $\{\Gamma^k\}$ satisfying (SS) in $U$. Assume that the metrics $g^k$ converge smoothly to a metric $g$, each $\Gamma^k$ is stable and minimal relative to the metric $g^k$, and $\sup \mathcal{H}^{n-1}(\Gamma^k) < \infty$. Then there are a subsequence of $\{\Gamma^k\}$ (not relabeled), a stable stationary varifold $V$ in $U$ (relative to the metric $g$), and a closed set $S$ of Hausdorff dimension at most $n-8$ such that
\begin{itemize}
\item[(a)] $V$ is a smooth embedded hypersurface in $U \setminus S$;
\item[(b)] $\Gamma^k \rightarrow V$ in the sense of varifolds in $U$;
\item[(c)] $\Gamma^k$ converges smoothly to $V$ on every $U' \subset\subset U \setminus S$.
\end{itemize}

\end{theorem}

\begin{remarq} The smooth convergence of the subsquence $\{\Gamma^k\}$ in part (c) is understood in the following sense: take an open set $U''\subset U'$ where the varifold $V$ is an integer multiple $N$ of  a smooth oriented surface $\Sigma$. Then, for $k$ sufficiently large, $\Gamma^k\cap U''$ is the union of $N$ disjoint surfaces $\Gamma_i^k$, $i=1,\dots,N$, that are normal graphs over $\Sigma$ of functions $f_i^k\in C^\infty(\Sigma)$. The convergence is smooth, in the sense that, for every $l\in\N$ and $\e>0$, $\|f_i^k\|_{C^l}<\e$, if $k$ is sufficiently large.
\end{remarq}

\begin{remarq} If $3\leq n\leq 7$ the closed set S from the theorem above is empty, in particular the limit $V$ is an embedded smooth hypersurface.
\end{remarq}

\

\

\section{Convergence of phase transitions}\label{convergence}

In this section we show that the regularity results for stable limit-interfaces from \cite{TonegawaWickramasekera}, can be extended to the case of bounded Morse index (Theorem \ref{morseindex}). 

\subsection{Assumptions}\label{assumptions} Let $M^n$ be a compact Riemannian manifold, $n\geq 3$, and $W\in C^3(\R)$. From now on we assume that we are in the following situation:
\begin{itemize}
\item [A.\ ] \textit{$W$ is a double-well potential}: $W\geq 0$, with exactly three critical points, two of which are non-degenerated minima at $\pm1$, with ${W(\pm1)=0}$ and $W''(\pm1)>0$, and  the third a local maximum $\gamma\in(-1,1)$.

\ \

\item [B.\ ] There are sequences $\e_k\to 0$ and $u_k\in C^3(M)$, such that $\sup_k E_{k}(u_k)<\infty$, $sup_{k}\|u_k \|_{L^\infty(M)}<\infty$ and $u_k$ is a critical point of the energy functional $$E_k(u)=\int_M \e_k \frac{|\nabla u|^2}{2} + \frac{W(u)}{\e_k},$$ i.e. it satisfies \begin{align}\label{equationk}-\e_k^2 \Delta u_k + W'(u_k)=0,\end{align} where $\Delta$ is the Laplace-Beltrami operator on $M$.

\ \
 
\item [C.\ ] There exists $m\in \N$, such that every $u_k$ has Morse index $m(u_k)\leq m$, i.e. the dimension of any subspace of $H^1(M)$ where the quadratic form $$E''_k(\phi,\phi)=\int_M \e_k |\nabla \phi|^2 + \frac{W''(u_k)}{\e_k}\phi^2$$ is negative definite, is at most $m$.
\end{itemize}

\subsection{Definition}\label{stable} Given an open set $U\subset M$, we say that $u_k\in C^3(M)$ is a \textit{stable critical point of $E_k$ in $U$}, if $u_k$ satisfies ($\ref{equationk}$) and $$\frac{d^2}{dt^2}E_k(u_k+t \phi)\bigg|_{t=0} = E''_k(\phi,\phi)=\int_M \e_k |\nabla \phi|^2 + \frac{W''(u_k)}{\e_k}\phi^2\geq 0$$ for every $\phi\in C^1(U)$.  

\subsection{Remark}\label{indexstable} An important immediate consequence of Assumption C, is that given any $(m+1)$-uple of disjoint open subsets of $M$, $u_k$ is a stable critical point of $E_k$ in at least one of them. To see this we can argue by contradiction. If there are functions $\phi_i \in C^1(M)$, $i=1,\dots,m+1$, with disjoint supports and such that $E_k ''(u_k)(\phi_i,\phi_i)<0$, they generate a $(m+1)$-dimensional subspace in which $E_k ''(u_k)(\cdot,\cdot)$ is negative definite, which contradicts $m(u_k)\leq m$.

\subsection{The associated varifolds} 

Given a sequence of critical points $u_k\in C^3(M)$ for the functional $E_k$, we associate to it a sequence of varifolds as in \cite{HutchinsonTonegawa}. Recall that if $\Sigma$ is a $(n-1)$-rectifiable subset of $M$,  $V_\Sigma$ denotes the varifold canonically induced by $\Sigma$.

Set $w_k=\Psi\circ u_k$, where $\Psi(t)=\int_0 ^t \sqrt{W(s)/2}ds$. We define $V_k$, the \textit{associated varifold to} $u_k$, by $$V_k(A)=\frac{1}{\sigma}\int_{-\infty}^\infty V_{\{w_k =t\}}(A)dt,$$ for every borel set $A$ and $\sigma=\int_{-1}^{1} \sqrt{W(s)/2}ds$. 

Notice that by the coarea formula $$\|V_k\|(A)=\frac{1}{\sigma}\int_{A}|\nabla w_k|=\frac{1}{\sigma}\int_A \sqrt{W(s)/2} \cdot |\nabla u_k|,$$ so one may interpret $V_k$ as a normalized averaging of the level sets of $u_k$.

\subsection{Optimal regularity}
If $V=V_{\theta\Sigma}$ is an integer rectifiable varifold, $\operatorname{reg} V$ denotes its \textit{regular set} (i.e. the points where $\Sigma$ is an embedded smooth hypersurface) and $\operatorname{sing}V$ denotes its \textit{singular set} (i.e. the complement of $\operatorname{reg} V$). 

Regularity problems deal with showing that under certain circumstances $\operatorname{sing} V$ is a small set. It is not true in general that an area minimizing varifold satisfies $\operatorname{sing} V=\emptyset$, but one can still show that it is a very small set. In this section we adopt the following notation.

\subsection{Definition}\label{optimalregularity} We say that a stationary integer rectifiable varifold $V=V_{\theta\Sigma}$ has \textit{optimal regularity}, if $\operatorname{sing} V$ has Hausdorff dimension at most $n-8$ (i.e. $\mathcal{H}^{n-8+\gamma}(\operatorname{sing} V)=0, \text{ for all } \gamma>0$). In particular  $\operatorname{sing} V$ is empty if $3\leq n\leq 7$ and $\operatorname{reg} V$ is an embedded minimal hypersurface. If in addition $\operatorname{reg}V$ is stable we will say that $V$ is \textit{stable with optimal regularity}.

\

For bounded open subset of $\R^n$, parts 1 and 2 of the following theorem are due to Hutchinson-Tonegawa \cite{HutchinsonTonegawa} and part 3 to Tonegawa-Wickramasekera \cite{TonegawaWickramasekera}. The same statements are true in closed manifolds (see Appendix B).

\subsection{Theorem}{\cite{HutchinsonTonegawa,Tonegawa,TonegawaWickramasekera}}\label{regularity} \ \textit{Suppose that Assumptions A and B from \ref{assumptions} hold. Then, after perhaps passing to a subsequence, we have:
\begin{enumerate}
\item \label{regularity1}The associated varifolds $V_{k}$ converge, in the varifold sense, to a stationary integral varifold $V$;
\item $\|V\|=\frac{1}{2 \sigma}\lim_{k\to\infty}E_k(u_k)$;
\item \label{regularity2} If in addition, the $u_k$ are stable critical points of $E_{k}$ on an open set $U$, then  $V\cap U$ is stable with optimal regularity.
\end{enumerate}
}

In item (3) they assume that $u_k$ are stable solutions of (\ref{eq:theequation}), i.e. $m(u_k)=0$. We extend the regularity result to the bounded Morse index case, i.e. Assumption C. This is the main result of this section and it implies Theorem A from the introduction.

\subsection{Theorem}\label{morseindex}\textit{Suppose that Assumptions A, B and C from \ref{assumptions} hold. Then, the varifold $V$ from parts (1) and (2) of Theorem \ref{regularity}, has optimal regularity.
}

\begin{proof}

By the remark after Definition \ref{stable}, given any $(m+1)$-uple of disjoint open subsets of $M$, each of the functions $u_k$ is stable in at least one of them. In particular there is a subsequence that is stable at least in one of them and part \ref{regularity2} of Theorem \ref{regularity} implies

\subsection*{Claim}\label{summarizing}\textit{Suppose that Assumptions A, B and C from \ref{assumptions} hold. Let $V$ be the varifold from parts (1) and (2) of Theorem \ref{regularity}. Then, given any $(m+1)$-uple of disjoint open subsets of $M$, at least in one of them $V$ is stable with optimal regularity.}

\

We expect the set $\operatorname{supp}V$ to have optimal regularity in all of $M$. As a first step in this direction, we show regularity in small annuli centered at an arbitrary point of $M$, following a similar argument as in the proof of Proposition $2.4$ from \cite{DelellisTasnady}. 

\subsection{Lemma}\label{annuli}\textit{There is a positive function $r:M\rightarrow \R^+$ such that in every annulus $An\in \mathcal{AN}(x,r(x))$, the varifold $V$ obtained in Corollary \ref{summarizing} is stable with optimal regularity.}

\begin{proof}
We argue by contradiction. Fix $x\in M$ and $0<\rho<\operatorname{Inj}(M)$, and suppose that for every $0<r<\rho$ there exists an annulus in $\mathcal{AN}(x,r)$ in which $V$ does not satisfies the conclusion of the lemma. Then by taking $r$ arbitrary small, we can find $m+1$ concentric, disjoint annuli centered at $x$ such that $V$ does not satisfies the conclusion of the lemma in any of them. This contradicts the claim above. 
\end{proof}

If $n\geq 8$, Lemma \ref{annuli} immediately implies that $V$ has optimal regularitiy, since we can \textit{hide} the potentially singular point $x$, inside the singular set of $\|V\|$. Then, we only need to explain why, if $3\leq n\leq 7$, the set $\operatorname{supp}\|V\|$ is also an embedded minimal hypersurface at the point $x$. This can be done  as in the proof of Lemma $5.2$ and Step $5$ of Proposition $2.8$ from \cite{DelellisTasnady}. We present here a sketch of the proof for convenience of the reader.

\subsection{Proposition} \textit{If $3\leq n\leq 7$, for every $x\in \operatorname{supp} V$, any tangent cone to $V$ at $x$ is an integer multiple of a hyperplane. Furthermore, $x$ is in the regular set of $V$, in particular $\operatorname{supp} V$ is a smooth embedded minimal surface.}

\begin{proof}
If $\rho$ is small enough, the varifold $V$ is a stable embedded minimal hypersurface (with integer multiplicity) $\Sigma$ in the punctured normal ball $B(x,\rho)\setminus\{x\}$. Let $T_\rho : \mathcal{B}_1(x) \to B(x,\rho)$ be the rescaled exponential map $T_\rho(z)=\exp_x(\rho z)$.

Given a sequence $\rho_k\to 0$, denote by ${\Sigma}_{k}$ the surface $T_{\rho_k}^{-1} (\Sigma)$. Then, by Theorem \ref{stablecompacity}, for any such a sequence, and every $\lambda>0$, there is a subsequence $\{\Sigma_{k_n}\}$ converging to a stable minimal hypersurface in the annulus $\mathcal{B}_{1-\lambda}(x)\setminus \overline{\mathcal{B}_{\lambda}(x)}$. In particular, since $\rho_k$ and $\lambda$ are arbitrary, any tangent cone to $V$ at $x$, is a stable minimal hypersurface in the punctured ball $\mathcal{B}_1 (x)\setminus\{0\}$, and by Simons' Theorem (see Theorem B.2 in \cite{SimonBook}), it must be an integer multiple of a hyperplane, because $n\leq 7$.

Choose $\rho_k=2^{-k}$. Then, given any positive constant $c_0$, for $k$ large enough, there is a plane $\pi_k$ such that $\Sigma_k \cap (\mathcal{B}_1(x)\setminus \overline{\mathcal{B}_{1/2}(x)})$ is the union of $m(k)$ disjoint graphs of Lipschitz functions over $\pi_k$, with Lipschitz constants smaller than $c_0$, counted with multiplicities $j_1(k),\dots,j_m(k)$, with $j_1+\cdots+j_m=N$.

We do not know a priori that the $\pi_k$ are all the same, but by comparing it with the tangent spaces of the Lipschitz graphs, it follows that the tilt between consecutive planes gets smaller as $k$ grows. In particular, in the annulus $An(x,2^{-k-3},2^{-k})$ the corresponding hypersurfaces of consecutive $k$'s must coincide, implying that they have the same number of components with the same multiplicities.

Doing this inductively, we find that $\Sigma\setminus\{x\}$ is the union of $m$ disjoint smooth embedded minimal hypersurfaces $\Gamma_1,\dots,\Gamma_m$, each homeomorphic to a disk minus a point and with multiplicities $j_1+\cdots+j_m=N$. Each tangent cone to $\Gamma_i$ is a hyperplane, and each $\Gamma_i$ is a minimal hypersurface with density 1.  It follows from Allard's regularity theorem that each $\Gamma_i$ is regular. Finally, since the $\Gamma_i$ are disjoint, $m>1$  would contradict the classical maximum principle. Thus $m=1$ and $x$ is a regular point for $\Sigma$.
\end{proof}

\end{proof}

\

\

\section{Min-max for phase transitions}\label{minmax.section}

In this section, we apply mountain-pass methods to construct a sequence $u_k\in C^3(M)$ of critical points of equation (\ref{equationk}) satisfying $-1\leq u\leq 1$ and with Morse index $m(u_k)\leq 1$. 

In what follows it will be convenient to modify the potential $W$ outside the set $[-1,1]$. More precisely, let $W^*\in C^3(\R)$ satisfy $W^*|_{[-1,1]}=W|_{[-1,1]}$, $W^*(x)>0$ for $|x|>1$ and constant on the set $\R\setminus [-2,2]$.  

Define the energy functional $E^*: H^{1}(M)\rightarrow \R$ as \begin{equation}\label{eq:energy} E^*(u)=\int_{M} \frac{|\nabla u|^2}{2}+W^*(u),\end{equation}

\subsection{Remark}\label{star} Notice that any $u$ that is a critical point of $E^*$ with $-1\leq u \leq 1$, is also a critical point of the functional $E(u)=\int_{M} \frac{|\nabla u|^2}{2}+W(u)$. Also, we have omitted any reference to the parameter $\e$, but all the results of this section apply to the functionals $E_k$ defined in the last section.

\ 

We can see right away that the functional $E^*$ is in a mountain-pass type situation. In fact, as a consequence of Assumption A, and the definition of $W^*$, the functions $\pm 1$ are the only global minimizers for $E^*$ in $H^1(M)$. On the other hand, the values of $E^*$ are bounded away from zero on the orthogonal complement of $\pm1$ in $H^1(M)$, i.e. the set of functions with zero average. More precisely:

\subsection{Lemma}\label{boundalpha}\textit{There is an $\alpha>0$ such that $E^*(u)\geq\alpha>0$ for every $u\in H^{1}(M)$ with $\int_M u =0$.
}

\begin{proof}
It is enough to show that there is a function $u$ with zero average with $E^*(u)=\alpha=\inf\{ E^*(u) : \int_M u=0 \}$. In fact, if $E^*(u)=\alpha=0$, then $u\equiv 1$ or $u\equiv -1$, which contradicts $\int_M u=0$.  

The existence of such a minimizer follows from a standard compactness argument. Take a sequence $u_n\in H^1(M)$, with $\int_M u_n =0$, and such that $E^*(u_n)\to\alpha$ . Since $\int_M u_n =0$, by the Poincare inequality, there is a constant $C>0$, such that $C\|u_n\|_{L^2(M)}\leq \|\nabla u_n\|_{L^2(M)}\leq E^*(u_n)$. In particular $u_n$ is a bounded sequence in $H^1(M)$, and by the Rellich-Kondrachov compactness theorem, there is an $u\in H^1(M)$, and a subsequence $\{u_{k}\} \subset \{u_{n}\}$ such that $u_k\to u$ weakly in $H^1(M)$ and strongly in $L^2(M)$. By the $L^2(M)$ convergence we must have $\int_M u = 0$. Then $E^*(u)\geq \alpha$. On the other hand, the functional $E^*$ is lower semicontinuos with respect to the weak convergence, this implies $E^*(u)\leq\lim E^*(u_k)=\alpha$. 
\end{proof}

This fact suggest that the subspace of functions with zero average, is a \textit{barrier} for the values of the energy between the points $\pm1$, making it plausible to obtain a critical point by min-max arguments. 

\subsection*{Setting for the min-max}\label{setting}
Let $\E$ be a $C^2$-functional on a Hilbert space $X=Y\oplus Z$, with $\dim(Y)=1$. In $Y$, identify the unit ball by $B_Y=\{y\in Y : \|y\|_{X} \leq1\}$ and the unit sphere by $S_Y:=\{y\in Y : \|y\|_{X} =1\}$. They consist of a line segment and two points, respectively. 

Assume that $Z$ is a barrier for the values of $\E$ between the points in $S_Y$, i.e. $$\alpha:=\inf \E|_{Z} >  \sup \E|_{S_Y}.$$
	
Let $\Gamma$ be the set of continuos paths with extrema in $S_Y$, i.e. $$\Gamma:=\{h:B_Y\rightarrow X \ | \ h \text{ is continous and } h|_{S_Y}=\operatorname{Id}_{S_Y} \}$$ and $c$ the \textit{min-max value} $$c:=\inf_{h\in\Gamma}\max_{t\in[0,1]}\E(h(t))\geq \alpha.$$ 

Denote by $K_c$ the set of critical points with energy $c$, i.e. $$K_c:=\{x\in X : \E(x)=c, \E'(x)=0\}.$$ Remember that for $x\in K_c$ its Morse index $m(x)$ is defined to be the index of the operator $\E''$, i.e. the maximal dimension of a subspace of $X$ where $\E''$ is negative definite.

A sequence $\{h_n\}$ in $\Gamma$ is called a \textit{minimizing sequence} if $$\max_{t\in[0,1]}\E(h_n(t))\to c \ \ \text{ as } \ n\to\infty.$$
	
Given a minimizing sequence $\{h_n\}$ we say that a sequence $\{x_n\}$ in $X$ is a \textit{min-max subsequence for $\{h_n\}$} if $$d(x_n,h_n(B_Y))_{X}\to 0 \ \ \text{and} \ \ \E(x_n)\to c,$$ as $n \to \infty$. 

Finally we say that $\E$ \textit{satisfies the Palais-Smale condition along $\{h_n\}$}, if every $\{x_n\}$ that is a min-max sequence for $\{h_n\}$ and satisfies $\E'(x_n)\to 0$, as $n\to\infty$, contains a convergent subsequence.

The following min-max theorem for functionals in a Hilbert space is of standard use in theory of semilinear elliptic partial differential equations. Its proof, in a much more general setting, can be found in the book \cite{Ghoussoub}. For our purposes it is enough to state a simplified version of Corollary 10.5 in \cite{Ghoussoub} adapted to our situation.

\subsection{Theorem}\label{min-max}\textit{Let $\E$ be the functional with the properties mentioned above. If $\E$ satisfies the Palais-Smale contidition along a minimizing sequence $\{h_n\}$ and if $\E''$ is Fredholm on $K_c$, then there exists $\{x_n\}$, a min-max subsubsequence for $\{h_n\}$, that converges to a critical point $x\in K_c$ with Morse index $m(x)\leq 1$.
} 

\ 

One advantage of Theorem \ref{min-max}, is that we only need to check the Palais-Smale condition along one minimizing sequence. It is possible to do this in the case of the functional $E^*$ defined in (\ref{eq:energy}). For this, and to verify the rest of its hypothesis, we rely on some well known properties of the functional $E^*$ that depend solely on the fact that the growth of the potential $W^*$ is controlled. We summarize them in the following 

\subsection{Proposition}\label{energyproperties}\textit{Let $E^*$ be the energy functional defined in} (\ref{eq:energy}), \textit{then:
\begin{enumerate} 
\item[i.] $E^*\in C^2(H^{1}(M))$ with derivatives 
\begin{subequations}
	\begin{align*}
	(E^*)'(u)(v) &=\int_M \nabla u \cdot \nabla v + (W^*)'(u) v  \\
	(E^*)''(u)(v,w) &=\int_M \nabla v \cdot \nabla w+ (W^*)''(u) v w. 
	\end{align*}
\end{subequations}
\item[ii.] $E^*$ satisfies the Palais-Smale condition for bounded sequences, i.e. if $u_n\in H^1(M)$, is a sequence such that $\|u_n\|_{H^1(M)}$ and $E^*(u_n)$ are bounded sequences, and $(E^*)'(u_n)\to0$, then it contains a convergent subsequence.
\end{enumerate}
}

\begin{proof}
The proof of (i) can be found on \cite{RabinowitzBook}, Appendix B, Propositions B.10 and B.34. To see (ii), we include here, for convenience of the reader, an sketch of \cite{ByeonRabinowitz}, Proposition 2.25.

Since the sequence $u_n$ is bounded in $H^1(M)$, by the Rellich-Kondrachov's compactness theorem, there is a $u\in H^1(M)$ and a subsequence $\{u_k\}\subset\{u_n\}$ such that $u_n$ is converging to $u$, weakly in $H^1(M)$ and strongly in $L^2(M)$.

We assert that $u$ is a critical point of $E^*$. In fact,
\begin{align*}
(E^*)'(u)(v) &=  \int \nabla u \cdot \nabla v + (W^*)'(u) v \\
& = \lim_n \int \nabla u_n \nabla v  + (W^*)'(u_n) v \\
& = \lim_n (E^*)'(u_n) (v) = 0, 
\end{align*}
since $u_n \to u$ in $L^2(M)$ and $(E^*)'(u_n)\to 0$ by hypothesis.

This implies that $(E^*)'(u_n)(u_n-u) - (E^*)'(u)(u_n-u)\to 0$, but on the other hand
$$(E^*)'(u_n)(u_n-u) - (E^*)'(u)(u_n-u)= \ \ $$
$$\ \  \int |\nabla(u_n-u)|^2 + \bigg[(W^*)'(u_n)-(W^*)'(u)\bigg](u_n-u),$$
and the second term on the right also goes to 0. 

In particular $$\int |\nabla (u_n-u)|^2 \to 0.$$

\end{proof}

As mentioned before, in order to apply the min-max theorem we just need one minimizing sequence such that the Palais-Smale condition holds along it. Since the functional $E^*$ satisfies the Palais-Smale condition on bounded sets of $H^1(M)$, it is enough to show that there exists a bounded minimizing sequence. 

Given any minimizing sequence $\{\tilde{h}_n\}_{n\in\N}$ we can obtain a bounded sequence by truncating it between the values $\pm1$. Define $h_n(t)=\min(\max(\tilde{h}_n(t),-1),1)$, for every $n\in\N$, observe that $-1\leq h_n(t)\leq 1$ and $E^*(h_n(t))\leq E^*(\tilde{h}_n(t))$. Then $\{h_n\}_{n\in\N}$ is also a minimizing sequence. Clearly $\|h_n(t)\|_{L^2(M)}$ and $\|\nabla h_n(t)\|_{L^2(M)}$ are bounded. Hence, the images of the paths $\{h_n\}_{n\in\N}$ are contained in a bounded subset of $H^1(M)$ and, by Proposition \ref{energyproperties}, the Palais-Smale condition is satisfied along this sequence.

Applying Theorem \ref{min-max} to this minimizing sequence, and by Remark \ref{star}, we obtain
\subsection{Proposition}\label{solution}\textit{There exists a function $u\in H^1(M)\cap K_c$, with ${-1\leq u\leq 1}$, and Morse index ${m(u)\leq 1}$. In particular $E(u)=c$ and, by standard elliptic regularity, $u\in C^3(M)$ is a classical solution to equation (\ref{eq:theequation}).
}

\

\ 

\

\section{A minimal hypersurface as a limit-interface}\label{limit-interface.section}

In this section we construct a minimal hypersurface as the limit-interface of a sequence of critical points for the functionals $$E_\e (u)=\int_{M} \e \frac{|\nabla u|^2}{2}+\frac{W(u)}{\e}.$$  The results from the last section apply to the functional $E_\e$, for each fixed $\e>0$. In particular, if we define \begin{equation}\label{minmaxval} c_\e:= \inf_{h\in\Gamma} \ \max_{t\in[0,1]}  \ E_\e(h(t)),\end{equation} we have

\subsection{Proposition}\label{solutionse}\textit{For every $\e>0$, there exist $u_\e \in C^3(M)$, a critical point for $E_\e$, with $E_\e(u_\e)=c_\e$, $-1\leq u_\e \leq1$, and Morse index $m(u_\e)\leq 1$.}

\

In order to apply Theorem \ref{morseindex} to a sequence of these $u_\e$ to obtain a non-trivial minimal hypersurface as a limit-interface, we need to verify that the energies $c_\e$ do not explode or vanish as $\e\to 0$. We have the following 

\subsection{Proposition}\label{bounds}
$$0<\liminf_{\e\to 0}c_\e \leq \limsup_{\e\to 0}c_\e<\infty.$$
\

The proof of Proposition \ref{bounds} is given in Section \ref{bounds.section} (lower bound) and Section \ref{upperbound.section} (upper bound). Here we present a sketch of the proof to motivate the computations done in those sections, but before that, lets state the main theorem of this section, that correspond with Theorem B of the introduction.

\subsection{Theorem}\label{phasetransitionminimal} \textit{In every $n$-dimensional compact manifold there exists an integral varifold $V$ such that
\begin{enumerate}
\item [(i)] $\|V\|=\frac{1}{2\sigma}\liminf c_\e$;
\item [(ii)] $V$ is stationary in $M$;
\item [(iii)] $\mathcal{H}^{n-8+\gamma}(\operatorname{sing}(V))=0,$ for every $\gamma>0$;
\item [(iv)] $\operatorname{reg}(V)$ is an embedded minimal hypersurface.
\end{enumerate}
}

\begin{proof}
By Proposition \ref{solutionse}, to every $\e>0$ we can associate a function $u_\e\in C^3(M)$, with $-1\leq u_\e \leq 1$, $m(u_\e \leq 1$, that is a critical point for the functional $E_\e$ at the level $E_\e(u_\e)=c_\e$. Any sequence $\{u_{\e_k}\}_{k=1}^\infty \subset \{u_\e\}_{\e>0}$ with energies converging to the value  $\liminf_{\e\to 0}c_\e>0$, will also have uniformly bounded energies by Proposition \ref{bounds}. We obtain the result after applying Theorem \ref{morseindex} to the sequence $u_{\e_k}$.
\end{proof}

\subsection{Sketch of the proof of Proposition \ref{bounds}}\label{sketch}

\begin{figure}\label{fig:iso} 
\includegraphics[width=220pt]{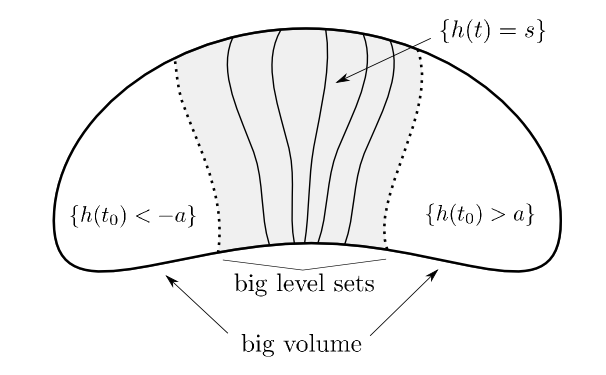}
\caption{Isoperimetric inequality}
\end{figure}

To see that $\liminf  c_\e$ is bounded away from zero, we use an isoperimetric inequality-type argument. The idea is to show that every path $h(t)$ joining $\pm1$ passes through a function with high energy. To see this choose $t_0$ such that $h(t_0)$ has zero average. On one hand, by the form of the potential $W$, the set $\{-a<h(t_0)<a\}$ has to be small. This implies that the set $\{h(t_0)\leq -a\}\cup\{h(t_0)\geq a\}$ is big. On the other hand, since the function $h(t_0)$ is bounded and has zero average, both sets will have to be big. By the isoperimetric inequality, the level sets $\{h(t_0)=s\}$ are big for $s\in[-a,a]$ (Figure 1). We obtain the lower bound for the energy after applying the coarea formula. This is done in Section \ref{bounds.section}.

To compute an upper bound for $\limsup c_\e$, we construct a path $h_\e(t)$ in $H^1(M)$ joining the constant functions $\pm 1$, with energy controlled independently of $\e$.

\begin{figure} 
\includegraphics[width=220pt]{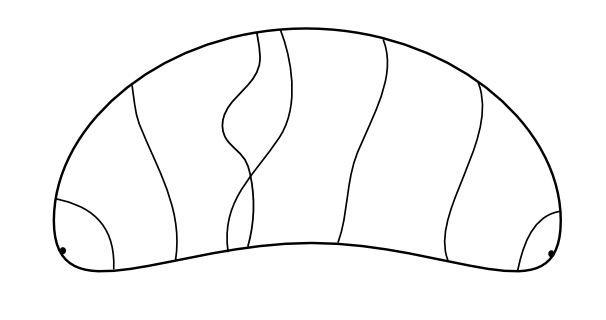}
\caption{Sweepout of $M$}
\end{figure}

Starting with a sweepout of $M$ by hypersurfaces $\{\Sigma_t\}$ coming from isotopic deformations of level sets of a Morse function (Figure 2),  we choose a small tubular neighborhood $N_t$ of $\Sigma_t$, and associate to it a function $h_\e(t)\in H^1(M)$ that is $+1$ in one component of $M\setminus N_t$ and $-1$ on the other (Figure 3), while in the tubular neighborhood the $h_\e(t)$ will have the profile of a 1-dimensional solution to the elliptic Allen-Cahn equation (Figure 4). 

\begin{figure} 
\includegraphics[width=220pt]{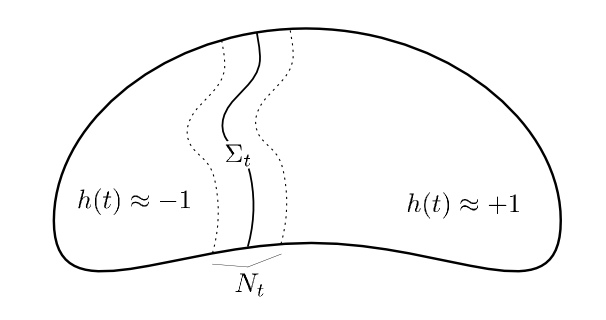}
\caption{The function $h(t)$ associated with the slice $\Sigma_t$}
\end{figure}

\begin{figure} 
\includegraphics[width=180pt]{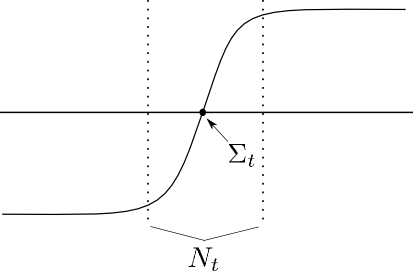}
\caption{Profile of $h(t)$ near $\Sigma_t$}
\end{figure}

Finally we show that the energies satisfy $$ E_\e(h_\e(t)) \to 2 \sigma \cdot \mathcal{H}^{n-1}(\Sigma_t),$$ as $\e\to 0$, uniformly on $t$. In particular we will have that $$\limsup c_\e \leq 2 \sigma \cdot \max_{t\in[0,1]}\mathcal{H}^{n-1}(\Sigma_t),$$ which is an upper bound indepedent of $\e$. This is done in Section \ref{upperbound.section}.


\

\

\section{Proof of Proposition \ref{bounds}: Lower bound}\label{bounds.section}

In this section we prove the lower bound in Proposition \ref{bounds}. The proof is an adaptation of arguments from \cite{ByeonRabinowitz} which is based on a local isoperimetric-type inequality due to De Giorgi.  

\subsection*{De Giorgi's Isoperimetric Inequality} The following lemma can be interpreted as a Sobolev space version of the isoperimetric inequality. A similar argument was used by De Giorgi in the proof of the regularity of solutions to elliptic PDEs. Roughly speaking, it states that functions in $H^{1}(M)$ cannot have jump singularities (see \cite{CaffarelliVasseur}).

\subsection{Lemma}\label{DeGiorgi}\textit{Let $u\in H^{1}(M)$ and suppose there are numbers $a < b$ such that ${\rm Vol}(\{u< a\})>\delta$ and ${\rm Vol}(\{b<u\})>\delta$, for some $\delta>0$. Then there is a positive constant $C=C(\delta,M)>0$, such that $$C (b-a)\ \leq {\rm Vol}(\{a\leq u\leq b\})^{1/2} \cdot \| \nabla u\|_{L^2(M)}$$}

\begin{proof}Given a compact manifold $M$, the function $\mathcal{I}:[0,\operatorname{Vol}(M)]\rightarrow\R$, defined by $$\mathcal{I}(t)=\inf \{\mathcal{H}^{n-1}(\partial \Omega) : \Omega\subset M \text{ and }\operatorname{Vol}(\Omega)=t \},$$ where $\Omega$ varies among all the sets of finite perimeter, is called the \textit{isoperimetric profile of} $M$. It is well known that $\mathcal{I}$ is continuous, vanishes on the extrema and is positive elsewhere.

Define $\Omega_t = \{u \leq t\}$, then for $t\in(a,b)$, we have ${\rm Vol}(\Omega_t)\in(\delta,{\rm Vol}(M)-\delta)$. By the continuity of $\mathcal{I}$ there exists a constant $C=C(\delta,M)>0$ such that $I(t)\geq C$ for such $t$. 
The set $\Omega_t$ has finite perimeter for almost every $t$ (see \cite{EvansGariepy}) and by the coarea formula 

\begin{align*}C(b-a) \leq & \int_a ^b \mathcal{H}^{n-1}(\partial \Omega_t) dt \\ = & \int_{\{a\leq u\leq b\}}|\nabla u| \\ \leq & {\rm Vol}(\{a\leq u\leq b\})^{1/2} \| \nabla u\|_{L^2(M)}\end{align*}

\end{proof}

By Lemma \ref{boundalpha} we know that $c_\e>0$, but \textit{a priori} we have no control over the behavior of $c_\e$ as $\e\to 0$. In what follows, we use the isoperimetric inequality from Lemma \ref{DeGiorgi} to guarantee that $ \liminf c_\e$ does not vanish.

We argue by contradiction. Suppose $\liminf c_{\e} =0$ and take a sequence $\e_k\to 0$, such that $c_{\e_k}\to 0$, as $k\to\infty$. For convenience we suppress any reference to the parameter $k$ along the proof. 

Fixed $\e>0$, choose a continuous path $h:[0,1]\to H^{1}(M)$, joining $\pm1$, with $-1\leq h(t) \leq 1$ and such that $$\max E_\e(h(t)) \leq c_\e +\e.$$ Select $t$ such that the function $u=h(t)$ has zero average, i.e. $\int_M u=0$. We assert that this function has high energy.

Fix $0<a<1$, by the form of the potential $W$, there is a constant $C_a>0$, depending only on $a$ and $W$, such that $W(u)\geq C_a$ in $\{-a\leq u\leq a\}$, then 
$$C_a {\rm Vol}(\{-a\leq u\leq a\}) \leq \int_{\{-a\leq u\leq a\}} W(u) \leq \e (c_\e +\e).$$

It follows from the last inequality, and $-1<u<1$, that $$0=\int_M u \ \ d\operatorname{Vol}(M) \leq -a\ {\rm Vol}(\{u<-a\}) \ + \ {\rm Vol}(\{u>a\}) \ + C_a^{-1} \e (c_\e +\e),$$ and $$\operatorname{Vol}(M)\leq \ {\rm Vol}(\{u<-a\}) \ + \ {\rm Vol}(\{u>a\}) \ + C_a^{-1} \e (c_\e +\e) .$$ 

Combining both we obtain $${\rm Vol}(\{u>a\}) \geq \frac{a}{2}\operatorname{Vol}(M)- C_a^{-1} \e (c_\e +\e).$$ Hence, if  $\e$ is small enough, ${\rm Vol}(\{u>a\})\geq \frac{a}{3}\operatorname{Vol}(M)$. Similarly it can be shown that ${\rm Vol}(\{u<-a\})\geq \frac{a}{3}\operatorname{Vol}(M)$.

Finally by Lemma \ref{DeGiorgi}, there is a constant $C=C(a,M)>0$ such that 
\begin{align*}
0<2aC \leq {\rm Vol}(\{-a\leq u \leq a\})^{1/2} \cdot \|\nabla u \|_{L^2(M)} \leq \sqrt{2C_a^{-1}} (c_\e +\e).
\end{align*}
 
This contradicts $c_\e\to0$. 
\begin{flushright}$\square$\end{flushright}


\section{Proof of Proposition \ref{bounds}: Upper bound}\label{upperbound.section}

\subsection{Sweepouts coming from level sets of Morse functions}\label{morseswo}

In this section we consider sweepouts of $M$ by hypersurfaces, generated by isotopic deformations of the level sets of Morse functions. More precisely, let $$\Lambda=\bigg\{ \{\Sigma_t \}_{t\in[0,1]} : \Sigma_t= \psi_t(f^{-1}(t),1)  \bigg\},$$
 where $\psi_t$ and $f$ vary on the set $C^\infty([0,1], \Is(M))$ and the set of all Morse functions taking values in $[0,1]$, respectively.

Given $\{\Sigma_t\}\in \Lambda$, define $$\mathcal{F}(\{\Sigma_t\})=\max \mathcal{H}^{n-1}(\Sigma_t),$$ and the \textit{width of} $\Lambda$ as $$m_0(\Lambda)=\inf  \mathcal{F}(\{\Sigma_t\}),$$ where the infimum is taken among all $\{\Sigma_t\}\in \Lambda$.

In this section we show $$\frac{1}{2\sigma} \limsup c_\e\leq m_0(\Lambda).$$ Roughly speaking, given any sweepout $\Sigma'_t$ as mentioned, we produce paths $h_t$ with energy controlled by $\mathcal{F}(\{\Sigma_t '\})$. We develop the technical details of the sketch presented in \ref{sketch}.

\subsection{Distance functions}
Let $\{\Sigma_t\}_{t\in[0,1]}\in \Lambda$. Lets call $d_{\Sigma_t}$ the \textit{signed} distance function from $\Sigma_t$, where we choose the sign of $d_{\Sigma_t}$ in such a way that, $d_{\Sigma_t}$ varies continuously and $d_{\Sigma_0}$, $d_{\Sigma_1}$ are nonnegative and nonpositive, respectively. 

It is well known that $d_{\Sigma_t} \in W^{1,\infty}(M)$ and satisfies the Eikonal equation $|\nabla d_{\Sigma_t}|=1$. In addition, in Section \ref{appendix.section}, we show that if the map $t\rightarrow \Sigma_t$ continuous in the Hausdorff distance, the functions $d_{\Sigma_t}$ vary continously in $H^{1}(M)$. We will make use of this fact in what follows.

\subsection{One dimensional Allen-Cahn equation}\label{1dim}
Let $\psi$ be the solution to the 1-dimensional ODE $$\begin{cases} \psi'(s)&=\sqrt{2W(\psi(s))} \\ \psi(0) &=\gamma,\end{cases}$$ where $\gamma\in(-1,1)$ is the only critical point of $W$ in this interval (see \ref{assumptions}). The following properties of $\psi$ are easy to check and we left their proof to the reader.

\begin{enumerate}
\item [(i)] $\psi$ solves the 1-dimensional elliptic Allen-Cahn equation;
\item [(ii)] $\psi:\R\to(-1,1)$ and is monotone increasing;
\item [(iii)] $\psi(s)\to\pm1$, as $s\to\pm\infty$;
\item [(iv)] $s\ W(\psi(s))\to0, \text{ as } s\to\pm \infty$;
\item [(v)] $2\sigma = \int_\R (\psi')^2 /2 + W(\psi).$
\end{enumerate}

\subsection{Functions associated to the sweepout}

Using $d_{\Sigma_t}$ and $\psi$, we can construct a path in $H^1(M)$ joining the functions $\pm1$, with energy concentrated in a small tubular neighborhood of every $\Sigma_t$. The normal profile of this function will be a scaling of $\psi$, in order to fit most of the energy inside the tubular neighborhood. 

For every $\Sigma=\Sigma_t$, $\delta>0$  and $\e>0$, define 
$$v_{\e,\delta}(\Sigma,x)=
\begin{cases} 
\psi\big(d_\Sigma(x)/\e\big) &\text{if } |d_\Sigma(x)| \leq \delta \\
\psi\big(\operatorname{sgn} d_\Sigma \cdot \delta/\e \big) &\text{if } |d_\Sigma(x)| > \delta 
\end{cases}.$$

For fixed $\delta$ and $\e$, the functions $g(t)(x)=g_t(x)=v_{\e,\delta}(\Sigma_t,x)$, for $t\in[0,1]$, form a continuous curve in $H^{1}(M)$, since $\Sigma_t$ is varying continuously with respect to the Hausdorff distance (see Section \ref{appendix.section}). 

Notice that for the extremal values $t=0,1$ the functions $g(0)$ and $g(1)$ are not constant. Since we want to construct a path joining the constant functions $\pm1$, we need to attach another deformation at the extremes. Lets do this for $g(0)$. The construction for $g(1)$ is analogous. 

$\Sigma_0$ consists of a finite number of points $P$ and $g(0)=v_{\e,\delta}(P,\cdot)\geq \gamma$  by the choose of the sign of $d_{\Sigma_0}$. For $t\in[0,1]$, define $f_t(x):=(1-t)+t \cdot v_{\e,\delta}(p,x)$. Then $f_t$ is a continuous deformation from $1$ to $v_{\e,\delta}(p,\cdot)$.

Similarly we can construct a path $\bar{f}_t$, joining $g(1)$ with the constant function $-1$. Combining the three, we have a continuous $h:[0,1]\to H^{1}(M)$ with $h(0)\equiv 1$ and $h(0)\equiv -1$.

\subsection{Controlling the energy of $h$}

Lets see first that \begin{equation}\label{path1}\max_{t\in[0,1]}E_\e(h(t)) = \max_{t\in[0,1]}E_\e(g(t)).\end{equation} In fact, we show that the path $f_t$ constructed above, have energy at most $E_\e(g(0))$ (similarly $\bar f_t$ have energy controlled by  $E_\e(g(1))$). This implies (\ref{path1}), since $h$ is just the juxtaposition of the paths $f_t, g(t)$ and $
\bar f_t$.

Notice that $|\nabla f_t|\leq |\nabla v_{\e,\delta}(P,x)| =|\nabla g(0)|$. Also, since $d_P(x)$ does not change sign (it is positive in this case) the values of $v_{\e,\delta}$ are all between $[\gamma,1)$, but the function $W$ is strictly decreasing on this interval, then we have $W(f_t(x))\leq W(v_{\e,\delta}(p,x)).$ 

\begin{align*}
E_\e(f_t) & = \int_M \e \frac{1}{2} |\nabla f_t|^2 +\frac{W(f_t)}{\e} \ d\operatorname{Vol}(x)\\
& \leq \int_M \e \frac{1}{2} |\nabla v_{\e,\delta}(p,x)|^2 +\frac{W(v_{\e,\delta}(p,x))}{\e} \ d\operatorname{Vol}(x) \\ & = E_\e(v_{\e,\delta}(p,x))
\end{align*}
The path $f_t$ joins the constant function $1$ with $g(0)$, with energy along $f_t$ at most $E_\e(g(0))$, similarly the energy of $\bar{f}_t$ is at most $E_{\e}(g(1))$. This proves (\ref{path1}). 

Then, to control the energy of $h$, we only need to deal with $$\max_{t\in[0,1]}E_\e(g(t)).$$

The energy of $v_{\e,\delta}(\Sigma,\cdot)$ is given by
$$E_\e(v_{\e,\delta}(\Sigma))=\int_M \frac{\e}{2} |\nabla v_{\e,\delta}(\Sigma,x)|^2 + \frac{1}{\e}W(v_{\e,\delta}(\Sigma,x)) d\operatorname{Vol}(x),$$ and we can estimate its value in the disjoint sets $\{|d_\Sigma|> \delta\}$ and $\{|d_\Sigma|\leq \delta\}$.

Notice that
$$\nabla v_{\e,\delta}(\Sigma,x)=
\begin{cases} 
1/\e \cdot \psi'(d_\Sigma(x)/\e) \cdot \nabla d_\Sigma(x) &\text{if } |d_\Sigma(x)| \leq \delta \\
0 &\text{if } |d_\Sigma(x)| > \delta 
\end{cases}.$$

The first integral is given by

\begin{multline}\label{eq:far}\int_{\{|d_\Sigma|> \delta\}}  \frac{\e}{2} |\nabla v_{\e,\delta}(\Sigma,x)|^2 + \frac{1}{\e}W(v_{\e,\delta}(\Sigma,x)) d\operatorname{Vol}(x)\leq \\ \operatorname{Vol}(M)\frac{1}{\e}\bigg(W(\psi(-\delta/\e))+W(\psi(\delta/\e))\bigg)\end{multline}

To compute the second integral we use the coarea formula and  $|\nabla d_\Sigma|=1$.

\begin{align*}\int_{\{|d_\Sigma|\leq \delta\}}  \frac{\e}{2} |\nabla v_{\e,\delta}(\Sigma,x)|^2 + \frac{1}{\e}W(v_{\e,\delta}(\Sigma,x)) d\operatorname{Vol}(x) \\
=\int_{-\delta}^{\delta}\frac{1}{\e}   \bigg[ \frac{\psi'(s/\e)^2}{2} + W(\psi(s/\e))\bigg] \cdot  \mathcal{H}^{n-1}(\{d_\Sigma = s\}) ds \\
= \int_{-\delta/\e}^{\delta/\e}   \bigg[ \psi'(s)^2/2 + W(\psi(s))\bigg] \cdot \mathcal{H}^{n-1}(\{d_\Sigma = \e s\}) ds 
\end{align*}

Remember that $\Sigma_t$ are isotopic deformations of level sets of a Morse function. Then, by the results from Section \ref{appendix.section}, we have that for every $\eta>0$, there exists a $\delta_0>0$ such that \begin{equation}\label{uniform level}\mathcal{H}^{n-1}(\{d_{\Sigma_t} = s\})\leq (1+\eta)\mathcal{H}^{n-1}(\Sigma_t)\end{equation} for every $|s|\leq \delta_0$ and every $t\in[0,1]$. This, and (v) from \ref{1dim}, imply  $$\int_{\{|d_\Sigma|\leq \delta\}}  \frac{\e}{2} |\nabla v_{\e,\delta}(\Sigma,x)|^2 + \frac{1}{\e}W(v_{\e,\delta}(\Sigma,x)) d\operatorname{Vol}(x) \leq 2\sigma(1+\eta)\mathcal{H}^{n-1}(\Sigma),$$ with $\eta\to 0$ as $\delta \to 0$.

Finally observe that property (iv) from \ref{1dim}, implies that for any $\delta>0$ fixed, the right side of (\ref{eq:far}) vanishes as $\e\to 0$, independently of $t$. 

Summarizing, we have that for every $\eta>0$, there is $\e_0$ such that $\e<\e_0$ implies that $$c_\e \leq 2\sigma(1+\eta)\mathcal{F}(\{\Sigma_t\}).$$ In particular, $$\frac{1}{2\sigma} \limsup c_\e \leq \mathcal{F}(\{\Sigma_t\}),$$ for $\{\Sigma_t\}\in \Lambda$ arbitrary.

\begin{flushright}$\square$\end{flushright}

\

\


\section{Comparison with Almgren-Pitts min-max theory}\label{almgrenpitts.section}

In this section we compare our results with Almgren-Pitts' min-max theory. We construct discrete sweepouts with width controlled by $\liminf c_\e$. It follows from this that the minimal hypersurface obtained in Theorem \ref{phasetransitionminimal}, has area at least the hypersurface produced by Almgren-Pitts. This gives us a better lower bound for $\liminf_\e c_\e$ than the one obtained in Section \ref{upperbound.section}.

\subsection*{Notation} We follow the notation from \cite{MarquesNevesWillmore,MarquesNevesInfinite,Montezuma}. The reader can consult these references for a more detailed account on Almgren-Pitts' methods and applications.

\medskip

\begin{tabular}{lll}

$I(1,k)$ && the cell complex of $I=[0,1]$ whose 0-cells and \\
&& 1-cells are $[0],\dots,[3^{-k}],\dots,[1-3^{-k}],[1]$ and \\
&& $[0,3^{-k}],\dots,[1-3^{-k},1]$, respectively.\\
${\bf I}_k(M)$ && the set of $k$-dimensional integral currents in $M$.\\
$\Z(M)$ && the subspace of ${\bf I}_{n-1}(M)$ of closed currents.\\
$\mathcal{F}$ && the flat seminorm on ${\bf I}_{k}(M)$\\
$\M$ && the mass seminorm on ${\bf I}_{k}(M)$\\
$$

\end{tabular}

\subsection{Discrete setting}    

Given a map $\phi: I(1,k_i)_0\to \Z(M)$, its \textit{fineness} $\f(\phi)$, defined as $$\max \  \{ \textbf{M}(\phi(x)-\phi(y): x, y \text{ adjacent vertices in } I(1,k_i)_0\},$$ is a discrete counterpart of the notion of continuity.

Instead of considering continuous maps from $I$ into $\Z(M)$, Almgren-Pitts theory is concerned with sequences of discrete maps into $\Z(M)$ with fineness tending to zero.

\subsection{Homotopy notions} Given $\phi_i : I(1,k_i)_0\to \Z(M)$, for $i=1,2$, we say that \textit{$\phi_1$ and $\phi_2$ are 1-homotopic in} $(\Z(M;\M),\{0\})$, with fineness $\delta$, if there exists $k\in \N$ and a map $$\psi: I(1,k)_0\times I(1,k)_0\to \Z(M)$$ such that 
\begin{enumerate}
\item $\f(\psi)<\delta$;
\item $\psi([i-1],x)=\phi_i({\bf n}(k,k_i)(x)), i=1,2,$ for every $x\in I(1,k)_0$;
\item $\psi(x,[0])=\psi(x,[1])=0$ for every $x\in I(1,k)_0.$
\end{enumerate}
\

(See \cite{MarquesNevesWillmore}, Section 7.1, for the definition of ${\bf n}$).

\

\subsection{Definition} An \textit{$$(1,\textbf{M})-\text{homotopy sequence of mappings into } (\mathcal{Z}_{n-1}(M;\textbf{M}),\{0\})$$}
is a sequence of maps $\{\phi_i\}_{i\in \mathbb{N}}$
\begin{equation*}
\phi_i : I(1,k_i)_0 \rightarrow \mathcal{Z}_{n-1}(M),
\end{equation*}
such that $\phi_i$ is $1$-homotopic to $\phi_{i+1}$ in $(\mathcal{Z}_{n-1}(M;\textbf{M}),\{0\})$ with fineness $\delta_i$ and
\begin{enumerate}
\item[(i)] $\lim_{i\rightarrow \infty} \delta_i =0$;
\item[(ii)] $\sup\{\textbf{M}(\phi_i(x)) : x \in \text{dmn}(\phi_i) \text{ and } i \in \mathbb{N}\}< \infty$.
\end{enumerate}

\

There is also a notion of homotopy between two $(1,\textbf{M})$-homotopy sequences of mappings into $(\mathcal{Z}_{n-1}(M;\textbf{M}),\{0\})$. 

\subsection{Definition} We say that $S^1 = \{\phi^1_i\}_{i\in \mathbb{N}}$ \textit{is homotopic with} $S^2 = \{\phi^2_i\}_{i\in \mathbb{N}}$ if $\phi^1_i$ is $1$-homotopic to $\phi^2_i$ with fineness $\delta_i$ and $\lim_{i\rightarrow \infty} \delta_i =0$.

\

This defines an equivalence relation on the set of $(1,\textbf{M})$-homotopy sequences of mappings into $(\mathcal{Z}_{n-1}(M;\textbf{M}),\{0\})$. An equivalence class is called a $(1,\textbf{M})$-homotopy class of mappings into $(\mathcal{Z}_{n-1}(M;\textbf{M}),\{0\})$. We denote $\pi^{\#}_1(\mathcal{Z}_{n-1}(M;\textbf{M}),\{0\})$ for the set of homotopy classes. 

\subsection{Width}\label{width}
Let $\Pi \in \pi^{\#}_1(\mathcal{Z}_{n-1}(M;\textbf{M}),\{0\})$ be a homotopy class and $S =\{\phi_i\}_{i\in \mathbb{N}} \in \Pi$. We define
\begin{equation}
\textbf{L}(S) = \limsup_{i\rightarrow \infty} \max\{\textbf{M}(\phi_i(x)) : x \in \text{dmn}(\phi_i)\}.
\end{equation}

The \textit{width of $\Pi$} is the minimum $\textbf{L}(S)$ among all $S \in \Pi$, 
\begin{equation}
\textbf{L}(\Pi) = \inf\{\textbf{L}(S) : S \in \Pi\}. 
\end{equation}

\

Finally, we present the min-max theorem from Almgren-Pitts (see \cite{Pitts,MarquesNevesWillmore,MarquesNevesInfinite,Montezuma}). 

\subsection{Theorem}\label{almgrenpitts}
\textit{If $\Pi \in \pi^{\#}_1(\mathcal{Z}_{n-1}(M;\textbf{M}),\{0\})$ is a non-trivial homotopy class, then $\textbf{L}(\Pi)>0$ and there exists an integral varifold $V$ such that
\begin{enumerate}
\item[(i)] $\|V\|(M) = \textbf{L}(\Pi)$;
\item[(ii)] $V$ is stationary in $M$;
\item[(iii)] $\mathcal{H}^{n-8+\gamma}(\operatorname{sing}(V))=0$, for every $\gamma>0$.
\item[(iv)] $\operatorname{reg}(V)$ is an embedded minimal hypersurface. 
\end{enumerate} 
}

\subsection{Almgren's Isomorphism}\label{isomorphism} To use the min-max theorem below we need to produce a non-trivial homotopy class.  In his Ph.D. thesis, Almgren constructed an isomorphism $$F^{\#}_M: \pi^{\#}_1(\mathcal{Z}_{n-1}(M;\textbf{M}),\{0\}) \to H_n(M).$$ Besides being one of the original motivation for applying min-max techniques on the set $\pi^{\#}_1(\mathcal{Z}_{n-1}(M;\textbf{M}),\{0\})$, this isomorphism serve also as a tool for showing that certain homotopy classes are non-trivial. 

Formally, the isomorphism is constructed in the following way. Given $$\phi: I(1,k)_0 \to \Z(M),$$ select $A_j\in {\bf I}_n(M)$ with least mass, such that $\partial A_j = \phi([(j+1)3^{-k}])-\phi([j3^{-k}])$. Then $$F^{\#}_M(\phi)=\bigg[\sum_{j=0}^{3^{k}-1}A_j\bigg].$$ Since $\phi([0])=\phi([1])=0$, the boundary of the sum above is zero, in particular $F^{\#}_M(\phi)$ is an element of $\Z(M)$.

Of course, as it is, the map we just described is not well defined. First, the $A_j$'s mentioned above might not exist nor be unique. Second, we want to define $F^{\#}_M$ for elements in $\pi^{\#}_1(\mathcal{Z}_{n-1}(M;\textbf{M}),\{0\})$, which classes are represented by sequences of discrete functions, rather thant by just one $\phi$. The following results of Almgren show that all this is possible (compare with Lemma 3.2 from \cite{MarquesNevesInfinite}).

\subsection{Lemma}\label{uniqueiso}\textit{There are constants $\nu_M>0$ and $\rho_M\geq 1$, such that if $S,T\in \Z(M)$, and $\mathcal{F}(S,T)\leq \nu_M$, then there is a unique (isoperimetric) choice $A_j\in {\bf I}_n(M)$ such that} $$\partial A_j= S-T \textit{\ \ \ and\ \ \ } \M(A_j) \leq \rho_M \mathcal{F}(S,T).$$

\

This shows that $F^{\#}_M(\phi)$ is well defined. The work of Almgren also shows that if $$\phi_i:I(1,k)_0 \to \Z(M) \text{ for $i=1,2,$}$$ are homotopic in the discrete sense with fineness in the flat topology less than $\nu_M$, we have $F^{\#}_M(\phi_1) = F^{\#}_M(\phi_2).$ Then, for $\Pi \in \pi^{\#}_1(\mathcal{Z}_{n-1}(M;\textbf{M}),\{0\})$ we can define $F^{\#}_M(\Pi)$ by taking any representant $\{\phi_i\}_{i\in \N}\in \Pi$ and defining $$F^{\#}_M(\Pi)=F^{\#}_M(\phi_i),$$ for any $\phi_i$ with fineness $\f(\phi_i)<\nu_M$.

\subsection{A non-trivial homotopy class}

In what follows we show that there is a non-trivial homotopy class $\Pi$ with width controlled by the min-max energy $c_\e$ (Proposition \ref{nontrivialclass}). Roughly speaking, $${\bf L}(\Pi)\leq \frac{1}{2\sigma}\liminf c_\e.$$

We do this in the following way. Let $h_t$ be a continuous path in $H^1(M)$ joining the constant functions $\pm 1$. First, we choose a finite number of level sets $\Sigma_i=h_{t_i}^{-1}(s_i)$, for some $s_i\in (-1,1)$. Each $\Sigma_i$ is selected in such a way that its area is controlled by the energy of $h_i$. Also, we are able to make $\Sigma_i$ and $\Sigma_{i+1}$ arbitrarily close with respect to the flat norm, depending on $\e$. This is done in Proposition \ref{discrete}.

To obtain a class in $\pi^{\#}_1(\mathcal{Z}_{n-1}(M;\textbf{M}),\{0\})$ we need to produce sequences that are also fine in the mass norm. This can be done by a result from \cite{Pitts}, that state that closeness in the flat norm $\mathcal{F}(\Sigma_i,\Sigma_{i+1})$ implies the existence of discrete deformations between $\Sigma_i$ and $\Sigma_{i+1}$, without increasing the mass and arbitrarily fine in the mass norm. This is the content of Proposition \ref{interpolation}.

Finally, to see that the sequence obtained belongs to a non-trivial homotopy class we show that its image by Almgren's isomorphism is not trivial.

\subsection*{Technical lemmas} The following two lemmas are important consequences of the uniform bound on the energies $c_\e$ and the continuity of the paths $h_t\in H^1(M)$ . The proof of the second lemma is an adaptation of some arguments from \cite{ByeonRabinowitz}.

\subsection{Lemma}\label{closeinflat}\textit{Let $0<\delta<1$. Given $\e>0$, let $h_t \in \Gamma(M)$, such that $E_\e(h_t)\leq c_\e+\e$, for all $t\in[0,1]$, then $${\rm Vol}( \{ |h_t|\leq 1-\delta\})\leq C_\delta^{-1} \e(c_\e+\e),$$ where $$C_\delta=\min W|_{\{|x|<1-\delta\}}>0,$$ is a constant depending only on $\delta$.
}

\begin{proof}
By the form of the potential $C_\delta$ is a positive constant. Integrating $W\circ h_t$ on the set $\{ |h_t|\leq 1-\delta\}$ we get $$C_\delta \cdot {\rm Vol}(\{ |h_t|\leq 1-\delta\}) \leq \int_M W( h_t(x)) d{\rm Vol}_x \leq \e (c_\e + \e).$$
\end{proof}

\subsection{Lemma}\label{volumes} \textit{Let $0<\delta<1$ and $\alpha\in(-1+\delta,1-\delta)$. Given $\e>0$, let $h_t \in \Gamma(M)$, with $E_\e(h_t)\leq c_\e+\e$, for all $t\in[0,1]$. Define $$\Omega_t=\{x\in M : h_t(x)>\alpha\}$$ then there exists $\rho>0$ such that $${\rm Vol}(\Omega_t\setminus \Omega_s)\leq 2{C_\delta}^{-1} \e(c_\e+\e),$$ whenever $|s-t|\leq \rho$, where $C_\delta$ is the constant  from Lemma \ref{closeinflat}, depending only on $\delta$.}

\

\begin{proof}
We argue by contradiction. Suppose there exists $\e>0$ such that for all $\rho>0$, there are $t=t(\rho)$ and $s=s(\rho)$ satisfying $|s-t|<\rho$ and $${\rm Vol}(\Omega_t\setminus \Omega_s)> {2}{C_\delta^{-1}} \e(c_\e+\e).$$

Since $\alpha\in(-1+\delta,1-\delta)$ and $$\Omega_t\setminus \Omega_s \ \subset\  \Omega_t \cap \big( \{|h_s|\leq 1-\delta \} \cup \{h_s < -1+\delta \} \big),$$ by Lemma \ref{closeinflat}, \begin{equation}\label{eq:vol} {\rm Vol}(\Omega_t \cap  \{h_s < -1+\delta \}) \geq C_\delta^{-1} \e(c_\e+\e).\end{equation}

Notice that $$\Omega_t \cap  \{h_s < -1+\delta \} \subset X_{s,t} = \{x\in M : |h_t(x) - h_s(x)| \geq \alpha +1-\delta \},$$ then (\ref{eq:vol}) implies \begin{equation}\label{eq:vol2}{\rm Vol}(X_{s,t}) \geq C_\delta^{-1} \e(c_\e+\e) \end{equation} for  $s$ and $t$ arbitrarily close. 

This contradicts the continuity of $t\mapsto h_t$ in $H^1(M)$ since $$ (\alpha +1-\delta)^2 {\rm Vol}(X_{s,t}) = \int_{X_{s,t}} |h_t(x) - h_s(x)|^2 \leq \|h_t - h_s \|^2_{H^1(M)}$$

\end{proof}

\subsection{Construction of the non-trivial homotopy class}

Define $\Psi(t)=\int_{s_0}^t \sqrt{W(t)/2}$, where $s_0$ is chosen so that $\Psi(\pm1)=\pm\sigma/2$. Given $h \in \Gamma(M)$ we define $\tilde h = \Psi \circ h$.  This is a normalized version of $h$, with values on the interval $[-\sigma/2,\sigma/2]$. Since $\Psi$ is strictly increasing, both functions have the same level sets.

Now we are ready to proof the following 

\subsection{Proposition}\label{discrete}
\textit{Given $\tilde \rho>0$ and $0<\tilde\delta<\sigma/2$, if $\e$ is small enough, there exists $k>0$, and a discrete sweepout $$\phi : I(1,k)_0 \to \mathcal{Z}_{n-1}(M)$$ such that 
\begin{enumerate} \item $\phi([0])=\phi([1])=0$
\item  $\mathcal{F}(\phi(a_j),\phi(a_{j+1}))\leq \tilde \rho$, for all $j=0,\dots,3^k-1,$ and $a_j=[j 3^{-k}]$.
\item $F^{\#}_M(\phi)\neq 0$.
\item ${\bf M}(\phi(a_j))< (c_\e + \e)/4 \tilde \delta$
\end{enumerate}
}

\begin{proof}
Given $\e>0$, let $h_t \in \Gamma(M)$, with $E_\e(h_t) < c_\e+\e$, for all $t\in[0,1]$. For every $t \in [0,1]$, choose $\tilde s(t)\in [-\tilde\delta,\tilde\delta]$ such that $\{\tilde h_t>\tilde{s}(t)\}$ is a set of finite perimeter and $\widetilde\Sigma_{t}=\partial \{\tilde h_t>s(t)\}$ satisfies $$2\tilde \delta \cdot \mathcal{H}^{n-1}(\widetilde\Sigma_t) \leq \int_{-\tilde\delta}^{\tilde\delta} \mathcal{H}^{n-1}(\partial \{\tilde h_t>s\}) ds \leq \int_M |\nabla \tilde h_t|.$$ 

It follows from the definition of $\tilde h$,  that $$|\nabla \tilde h_t|=|\nabla h_t|\sqrt{W(h_t)/2}\leq \frac{1}{2}\bigg(\e\frac{|\nabla h_t|^2}{2}+\frac{W(h_t)}{\e}\bigg).$$ In particular, \begin{equation}\label{swo1}\mathcal{H}^{n-1}(\widetilde\Sigma_t) < (c_\e+\e)/4\tilde\delta.\end{equation}

Notice that by the way $\tilde h$ was constructed, there are $\delta>0$ (depending only on $\tilde \delta$) and $s(t) \in (-1+\delta,1-\delta)$, such that $\widetilde\Sigma_{t}=\partial \{h_t>s(t)\}$. Let $\rho>0$ be the constant given by Lemma \ref{volumes}, and take $k\in\N$ such that $3^{-k}<\rho$. We assert that $\phi([j3^{-k}])=\widetilde \Sigma_{j3^{-k}}$ is the discrete sweepout we want.

(1) and (4) follows from $h_0\equiv -1$, $h_1\equiv 1$ and (\ref{swo1}), respectively. To see (2), choose $\alpha\in (-1+\delta,1-\delta)$ so that $\Omega_t=\{h_t > \alpha\}$ is a set of finite perimeter, for all $t\in \mathbb{Q}$, and define $\Sigma_{j3^{-k}}=\partial \Omega_{j3^{-k}}$, for all $j=0,1,\dots,3^k$.

As currents, $$\Sigma_t - \Sigma_s=\partial A({s,t}) \text{\ \  and}$$ $$\widetilde\Sigma_{j3^{-k}} - \Sigma_{j3^{-k}}=\partial B_j,$$ where $A({s,t})=\Omega_t-\Omega_s$ and $B_j=\{h_{j3^{-k}}>s(j3^{-k})\}-\Omega_{j3^{-k}}$. 

Also, for open sets $U$ and $V$ considered as currents, we have $\operatorname{supp}( U-V) \subset (U\setminus V) \cup (V\setminus U).$ Then, from Lemmas \ref{volumes} and \ref{closeinflat} it follows $$\mathcal{F}(\Sigma_{j3^{-k}},\Sigma_{(j+1)3^{-k}}) \leq 4 C_\delta^{-1}\e(c_\e+\e),\ \ \ \text{for $j=0,1,\cdots,3^k-1$ \ and}$$ $$\mathcal{F}(\widetilde\Sigma_{j3^{-k}},\Sigma_{j3^{-k}})\leq C_\delta^{-1}\e(c_\e+\e), \ \ \text{for all $j=0,1,\cdots,3^k$}.$$ 

Then $$\mathcal{F}(\widetilde\Sigma_{j3^{-k}},\widetilde \Sigma_{(j+1)3^{-k}})\leq 6 C_\delta^{-1}\e(c_\e+\e),$$ for all $j=0,1,\cdots,3^k-1$ and (2) follows choosing $\e$ sufficiently small.

Finally, to see (3) define $C_j=B_{(j+1)3^{-k}}+A((j+1)3^{-k},j3^{-k})-B_{j3^{-k}}$. Notice that if $\e$ is small enough, Lemma \ref{uniqueiso} implies, that $C_j$ is the only (small) current satisfying $\widetilde\Sigma_{(j+1)3^{-k}} -\widetilde\Sigma_{j3^{-k}}=\partial C_j$. Then $F_M^{\#}(\phi)$ is well defined by the formula $$F_M^{\#}(\phi)= \bigg[ \sum_{j=0}^{3^k-1} C_j \bigg]=\bigg[ \sum_{j=0}^{3^k-1} A_j \bigg]=\big[\Omega_{1}\big]-\big[\emptyset\big]=\big[M\big].$$
\end{proof}

Proposition \ref{discrete} provides a discrete map of currents, arbitrarily fine in the \textit{flat} norm, with controlled mass and non-trivial image under Almgren's isomorphism. However, to apply Almgren-Pitts min-max technique, we need to produce discrete maps that are fine in the \textit{mass} norm. This is the discrete analogue of a situation that is common on recent applications of the min-max technique, in which sweepouts continuous with respect to the flat norm arise naturally (see \cite{MarquesNevesWillmore,MarquesNevesInfinite,Rosenberg}).

Unfortunately, an important technical difficulty appears when trying to pass from the flat to the mass, due to the phenomenon of \textit{concentration of mass}. The problem is that a limit of currents can be quite different from a limit of varifolds, i.e. if $S_i\to S$ is a sequence of currents converging in the flat norm, such that the induced varifolds $|S_i|\to V$ are converging in the weak topology, it is not true in general that $|S|=V$. 

In \cite{MarquesNevesWillmore,MarquesNevesInfinite}, the notion of sweepouts with no concentration of mass was introduced to deal with this problem in the multiparameter min-max. Our situation if different for two reasons. On one hand, we are dealing with discrete maps rather than continuous sweepouts, and on the other hand we are only interested in the 1-parameter min-max. The results we need follow almost immediately from the work of Pitts \cite{Pitts}. 

Concerning the general case, in a recent work \cite{Zhou2015}, Xin Zhou showed that in fact it is not necessary to assume the no concentration of mass condition.

\subsection{Technical Lemma: No concentration of mass}\label{noconcentration} Lemma 3.7 from \cite{Pitts} allow us to rule out the case of concentration of mass, we state it here adapted to our context.

Let $T, T_1,T_2,\dots$ be elements in $\Z(M)$ and $V\in \mathcal{V}_{n-1}(M)$, such that \begin{enumerate}
\item $T_i\to T$ in the flat norm, and 
\item $|T_i|\to V$ as varifolds.
\end{enumerate} 
Then, for every $\delta>0$, there exists a sequence $S_1,S_2,\dots$ in $\Z(M)$ such that \begin{enumerate}
\item $S_i\to T$ in the flat norm, and 
\item $|S_i|\to |T|$ as varifolds,
\end{enumerate} 
and for every $i\in \N$, there is a finite sequence $R_0,\dots,R_m \in \Z(M)$,  such that 
\begin{itemize}
\item $R_0=T_i$ and $R_m=S_i$
\item $\M(R_i)\leq \M(T_i)+\delta, j=1,2,\dots,m$,
\item $\sup_j\M(R_j-R_{j-1})\leq \delta$.
\end{itemize}

\subsection{From the flat norm to the mass norm}\label{interpolation}

Lemma 3.8 from \cite{Pitts} states that currents close enough in the flat norm can be deformed one into another by means of a finite sequence of deformations arbitrarily fine in the mass norm. We state it here adapted to our context with some modifications from \cite{MarquesNevesWillmore}, Lemma 13.4. Define ${\textbf B}_s^{\mathcal{F}}(T)=\{ \tilde S\in \Z(M) : \mathcal{F}(T,\tilde S)\leq s \}$.

\subsection{Lemma}\label{bolaflat}\textit{Given $T\in \Z(M)$, $\delta>0$ and $L>0$, there exists $\nu=\nu(T,\delta,L)>0$ for which the following holds: Given $0<s<\nu$ and $S\in {\textbf B}_s^{\mathcal{F}}(T)\cup \{ \tilde S\in \Z(M): \M(\tilde S)\leq L\}$ then, for some $k\in \N$, there exists $$\phi:I(1,k)_0\to {\textbf B}_s^{\mathcal{F}}(T)$$ with 
\begin{enumerate}
\item $\phi([0])=S$ and $\phi([1])=T$,
\item ${\bf f}(\phi)\leq \delta$,
\item $\sup {\bf M}(\phi)\leq L+\delta$,
\end{enumerate}
}

\begin{proof}
Parts (1)-(3) follow from \ref{noconcentration} exactly as in \cite{Pitts}, Lemma 3.8. That $\phi([x])\in {\textbf B}_s^{\mathcal{F}}(T)$, for $x\in I(1,k)_0$, follows from the observation made just before formula (78) in the proof of Lemma 13.4 in \cite{MarquesNevesWillmore}.
\end{proof}

The fact that $\phi \in {\textbf B}_s^{\mathcal{F}}(T)$ is essential to guarantee that the image of the map by Almgren's isomorphism remains the same after refining it in the mass norm. In fact, notice that choosing $s$ in Lemma \ref{bolaflat} small enough, there exist $Q,Q_1,\dots,Q_k$, unique elements of ${\textbf I}_n(M)$, given by Lemma \ref{uniqueiso}, and satisying
\begin{itemize}
\item $\partial Q=T-S$
\item $\partial Q_i=\phi([i\cdot 3^k])-\phi([(i-1)\cdot 3^k])$, for all $i=1,\dots,k$.
\end{itemize}

If $s$ is small enough, we can also guarantee that $Q=Q_1+\cdots+Q_k$. In fact, let $\widetilde Q_i\in {\textbf I}_n(M)$ be the unique isoperimetric choice such that $\partial \widetilde Q_i = \phi([i\cdot 3^k])-S$. By definition $\widetilde Q_1=Q_1$. We also have that $\partial(Q_1+Q_2)=\phi([2\cdot 3^k])-S$, but $\M(Q_1+Q_2)\leq 2s$, and if $s$ is small enough Lemma \ref{uniqueiso} gives us $Q_1+Q_2=\widetilde Q_2$. Proceeding inductively and noticing that $\widetilde Q_k=Q$ we conclude that $$Q=Q_1+\cdots+Q_k.$$

Covering the space of bounded cycles with finite balls of this type $${\textbf B}_{s_i}^{\mathcal{F}}(T_i)\cup \{ \tilde S\in \Z(M): \M(\tilde S)\leq L\}$$ and arguing as in Lemma 3.8 of \cite{Pitts}, we conclude

\subsection{Lemma}\label{discreteinterpolation}\textit{
Fix $L>0$. There exists $\nu>0$ such that if $S,T\in\Z(M)\cup \{ \tilde S\in \Z(M): \M(\tilde S)\leq L\}$ satisfy $\mathcal{F}(S,T)<\nu$ there exists $k\in \N$ and $$\phi:I(1,k)_0\to\Z(M)$$ with 
\begin{enumerate}
\item $\phi([0])=S$ and $\phi([1])=T$,
\item ${\bf f}(\phi)\leq \delta$,
\item $\sup {\bf M}(\phi)\leq L+\delta$.
\end{enumerate}
Additionally if $Q,Q_1,\dots,Q_k$, are the unique elements of ${\textbf I}_n(M)$, given by Lemma \ref{uniqueiso}, and satisying
\begin{itemize}
\item $\partial Q=T-S$
\item $\partial Q_i=\phi([i\cdot 3^k])-\phi([(i-1)\cdot 3^k])$, for all $i\in \N$.
\end{itemize}
we have that $Q=Q_1+\cdots +Q_k$.
}

\

Combining the results from this section we obtain 

\subsection{Corollary}\label{finally}\textit{ Given $\tilde\delta<\sigma/2$, if $\e$ is small enough, there is a non-trivial homotopy class  $\Pi \in  \pi^{\#}_1(\mathcal{Z}_{n-1}(M;\textbf{M}),\{0\})$ such that $$ 0<{\bf L}(\Pi)\leq (c_\e+2\e)/4\tilde\delta.$$ 
}

\begin{proof}
From Proposition \ref{discrete}, for $\e$ small enough, we can find a map ${\phi:I(1,k)_0\to \Z(M)}$ with $$\sup \M(\phi)\leq (c_\e+\e)/4\tilde\delta \text{ \ \ and \ \ } \mathcal{F}(\phi(a_j),\phi(a_{j+1}))\leq \tilde \rho,$$ for all $j=0,\dots,3^k-1,$ where $a_j=[j 3^{-k}]$. 

Define $\tilde \phi:I(1,\tilde{k})_0\to \Z(M)$, as the refinement of $\phi$ given by applying Lemma \ref{discreteinterpolation} to every pair of adjacent vertices of $I(1,k)_0$, with $\sup \M(\tilde \phi)\leq (c_\e+2\e)/4\tilde\delta$. Choosing $\f(\tilde \phi)$ arbitrarily small there exists $\Phi$, the \textit{Almgren's extension} of $\tilde \phi$ (see Theorem 3.10 of \cite{MarquesNevesInfinite}). By the arguments in \ref{isomorphism} and the last part of Lemma \ref{discreteinterpolation}, we must have $$F_M(\Phi)=F_M^{\#}(\tilde \phi)=F_M^{\#}(\phi).$$ Then, $F_M(\Phi)=[M]$ by Proposition \ref{discrete}. 

\end{proof}

Now we are ready to prove Proposition D from the Introduction

\subsection{Proposition}\label{nontrivialclass} {\it There is a non-trivial class $\Pi$ such that $$0<{\bf L}(\Pi)\leq \frac{1}{2\sigma}\liminf c_\e.$$ In particular, if $V$ is the limit-interface of Theorem \ref{phasetransitionminimal}, and $V_{AP}$ is the stationary varifold obtained after applying Almgren-Pitt's min-max to the class $\Pi$ (Theorem \ref{almgrenpitts}), then $$\|V_{AP} \| \leq \|V\|.$$}

\begin{proof}
For $\tilde\delta<\sigma/2$  fixed, Corollary \ref{finally} holds for every $\e$ small enough, in particular, for the class $\Pi=(F_M^{\#})^{-1}([M])$, we have ${\bf L}(\Pi) \leq \frac{1}{4\tilde\delta}\lim\inf c_\e,$ and taking the limit as $\tilde \delta \to\sigma/2$, $${\bf L}(\Pi)\leq \frac{1}{2\sigma}\liminf c_\e.$$
\end{proof}

\subsection{Comparison with Almgren-Pitts' min-max theory}

In \cite{ColdingDelellis}, a continuous refinement of the Almgren-Pitts min-max theory is presented for $n=3$. In this context the min-max procedure is applied to sweepouts of $M$ by surfaces, coming from isotopic deformation of level sets of Morse functions.

We considered sweepouts of this kind in Section \ref{upperbound.section}. The results we presented, imply that $\frac{1}{2\sigma}\limsup c_\e$ is at most the area of the surface obtained in \cite{ColdingDelellis}. Combining this with the results from this section obtain Corollary E of the Introduction, i.e. for $n=3$, $$\|V_{AP}\|  \leq \| V\|  \leq \|V_{cont}\|,$$ where $V_{AP}, V$ and $V_{cont}$, are the varifolds obtained applying the Almgren-Pitts Theorem \ref{almgrenpitts} (to the class $\Pi$ of Proposition \ref{nontrivialclass}), the phase transitions Theorem \ref{phasetransitionminimal} and the continuous version of Almgren-Pitts from \cite{ColdingDelellis} (to the saturated family generated by the level sets of a Morse function), respectively. Then, in some sense, our phase-transitions construction of minimal surfaces, lies in between the original construction of Almgren-Pitts \cite{Pitts} and the continuous refinement presented in Colding-De Lellis \cite{ColdingDelellis}. 

When $n>3$, a similar refinement of the work of Almgren-Pitts is presented by De Lellis-Tasnady \cite{DelellisTasnady}. Nonetheless, the sweepouts they considered have much more singularities than the ones in \cite{ColdingDelellis}. This flexibility is essential to obtain the regularity results. The methods from Section \ref{upperbound.section} cannot be applied directly to these general class of sweepouts, so in this case we only have the inequality $$\|V_{AP}\|  \leq \| V\|.$$ 

Nonetheless, notice that in Section \ref{upperbound.section}, we use the fact that the sweepout comes from the level sets of Morse functions, only to obtain the uniform inequality (\ref{uniform level}) in the last part of the argument. The rest of the construction relies solely on the fact that the slices of the sweepout vary continuously with respect to the Hausdorff distance. We believe that the Morse function condition can be removed.


\section{Appendix A: Distance functions}\label{appendix.section}

In the first part of this section we prove some regularity properties of the distance function $d_K$ from a compact set $K$. We are specially interested in the behavior of $d_K$ when the set $K$ is moving continuously with respect to the Hausdorff distance. 

In the rest of the section we indicate how to obtain some estimates for the area of hypersurfaces parallel to a given one. Our objective is to show that the area of hypersurfaces parallel to the slices of a sweepout, can be chosen arbitrary close to the area of the slices in an uniform way along the sweepout. 

\subsection*{Regularity of the distance function}

Let $M$ be a compact Riemannian manifold and  $K$ a non-empty compact subset of $M$. Define $$d_K(x):=d(x,K):=\inf \{d(x,y):y\in K\},$$ where $d$ is the distance on $M$.

Also, given $K_1$ and $K_2$, compacts subsets of $M$, we define the \textit{Hausdorff distance} by $$d_H(K_1,K_2):=\max\{ \max_{x\in K_1}d(x,K_2) , \max_{x\in K_2}d(x,K_1) \}.$$ It is well known that $d_H$ is a metric on the set of all non-empty compact subsets of $M$, and clearly, $d_H(\{x\},K)=d_K(x)$.

Of course, $d_K$ is \textit{never} a smooth function in $M$. Nonetheless it possesses some  regularity properties that can be useful for applications. The first of the propositions below is a well known fact, but we present its proof since some parts of it are needed in the proof of the second proposition.

\subsection{Proposition}\label{eikonal} \textit{The function $d_K$ is differentiable almost everywhere on $M$ and satisfies the Eikonal equation $\|\nabla d_K\|=1$.}

\subsection{Proposition}\label{hausdorff} \textit{If $\{K_n\}_n$ is a sequence of non-empty compact subsets of $M$, converging in the Hausdorff distance to a compact set $K$, then, for every $1\leq p <\infty$, the functions $d_{K_n}$ converge to $d_K$, in $W^{1,p}(M)$, as $n\to\infty$.}

\subsection*{Proof of Proposition \ref{eikonal}} By the triangle inequality for $d_H$ we have $d_K(x)\leq d_K(y)+d(x,y)$, this implies that $|d_K(x)-d_K(y)|\leq d(x,y),$ so $d_K$ is a Lipschitz function (on $M$ with the metric $d$, but also in any chart as a function of $\R^n$). By Rademacher's theorem it is differentiable almost everywhere, and in fact $d_K\in W^{1,\infty}(M)$ with $\|\nabla d_K \|\leq 1$.

It is left to see that $\|\nabla d_K \|=1$, a.e. For any $x\in M$, there is at least one minimizing geodesic $\gamma : [0,d_k(x)]\rightarrow M$ joining $x$ and $K$. We prove the following about the points where $d_K$ is regular.

\subsection*{Claim} \textit{If $x$ is a regular point for $d_K$, the geodesic $\gamma$ mentioned in the last paragraph, is unique.}

\

If $\gamma : [0,d_k(x)]\rightarrow M$ is a minimizing geodesic joining $x$ and $y$, then $d_K(\gamma(s))=d_K(x)-s$ for every $s\in[0,d_K(x)]$, otherwise we would be able to find a shorter path joining $x$ to $K$, contradicting the minimality of $\gamma$. In particular $$\nabla d_K(x) \cdot \gamma'(0)=(d/ds) d_K(\gamma(s))\big|_{s=0}=-1,$$  since we already know that $\|\nabla d_K(x)\|\leq 1$ and $\|\gamma'(0)\|=1$, this can only happen if $\nabla d_K(x)=-\gamma'(0)$. In particular $\|\nabla d_K(x)\|=1$, and $\gamma$ is unique, proving the claim and the proposition. \begin{flushright}$\square$\end{flushright}

\subsection*{Proof of Proposition \ref{hausdorff}}
By the triangle inequality for $d_H$ we have $d(x,K)\leq d(x,K_n) + d_H(K_n,K)$. Interchanging the roles of $K_n$ and $K$, we obtain $$\sup |d(x,K) - d(x,{K_n})| \leq d_H(K_n,K).$$ In particular $d(x,K_n)\to d(x,K)$ in $L^\infty(M)$ (and also in $L^p(M)$, for $1\leq p < \infty$, since $M$ is compact). 

To see that $\nabla d_{K_n}$ converge to $\nabla d_{K}$ in $L^p(M)$, i.e. $\|\nabla d_{K_n}-\nabla d_{K}\|\to 0$ in $L^p(M)$, for $1\leq p < \infty$, observe that since we are in a compact manifold, and the norms $\|\nabla d_{K_n}\|$ and $\|\nabla d_{K}\|$ are bounded by 1, point-wise convergence a.e. will imply $L^p$-convergence, by the $L^p$-Lebesgue dominated convergence theorem. Then it is enough to proof that $\nabla d_{K_n}\to\nabla d_{K}$ a.e.

Proposition \ref{eikonal} implies that, for almost every $x\in M$, the function $d(\cdot,K_n)$ is differentiable in $x$, $\forall n\in\N$. By the Claim above, for such an $x$, and for every $n\in\N$ there is an unique geodesic $\gamma_n :[0,d_{K_n}(x)]\rightarrow M$ realizing the distance from $x$ to $K_n$. Lets call $\gamma$ the geodesic associated in the same way to the limit set $K$. 

We must have that $\gamma_n\to \gamma$, as $n\to\infty$. If not there would be another geodesic realizing the distance from $x$ to $K$, contradicting the uniqueness of $\gamma$. This implies that $\gamma_n'(0)\to \gamma'(0)$.

\begin{flushright}$\square$\end{flushright}

\subsection{Parallel hypersurfaces}\label{parallelhypersurfaces}

Let $M$ be a complete $n$-dimensional Riemannian manifold, and $\nu:N\to M$ an isometric embedding of a closed $(n-1)$-dimensional orientable manifold $N$. Let $n:N\to TM$ denote a choose of a normal vector field over $N$, i.e. $n(x)\in TM_{\nu(x)}$ and $n(x)\perp TN_{\nu(x)}\subset TM_{\nu(x)}$.

For such an $n$ we associate a \textit{normal exponential map} $\exp_n: N\times \R \to M$ given by $\exp_n(x,t)=\exp_{\nu(x)}(tn(x))$. Then $t\to \exp_n(x,t)$ is the geodesic emanating from $\nu(x)$ with velocity $n(x)$.

We are interested in estimating by above the size of the level sets $\{d_{\nu(N)}=\delta\}$, where $d_{\nu(N)}$ is a signed distance function to $\nu(N)$, positive in the direction of $n$. Note that $\{d_{\nu(N)}=\delta\}\subset \exp_n(N\times\{\delta\})$, so it is enough to give an upper bound for the area $$\int_N | \operatorname{Jac}_x (\exp_n)(x,\delta) | d\operatorname{Vol_N}(x)$$ where $\operatorname{Jac}_x (\exp_n)(x,\delta)$ denotes the Jacobian determinant of the map $\exp_n(\cdot,\delta):N\to M$ that sends $x\to \exp_n(x,\delta)$.

It is possible to give such estimates in terms of curvature bounds for $M$ and $N$. From now on we will assume that there are constants $k\geq 0$ and $\lambda \geq 0$, such that $$ -k< K_{min}\leq K_{\max}<k$$ and  $$|\langle S_{n(x)}(v),v\rangle|\leq \lambda \langle v, v\rangle$$ where $K_{\min}$ and $K_{\max}$ are the minimum and maximum  of the scalar curvatures in $M$, respectively, and $S_{n(x)}$ is the shape operator of the surface $N$ at $x$, associated to the normal vector $n(x)$.

The volume element of the map $\exp_n(x,t)$ is controlled by the norm of Jacobi vector fields along the geodesic $\exp_{x}{(t n(x))}$. The classical Rauch's comparison theorem can be extended to this setting with some modifications. In this case the Jacobi vector fields along the normal geodesic that one should consider does not vanish at initial time, as a consequence, its growth is not only controlled by the ambient sectional curvatures $K_{\min}$ and $K_{\max}$, but also by the initial conditions imposed by the shape operator $S_{n(x)}$ (see \cite{HeintzeKarcher,Warner}). The following proposition is an easy consequence of Corollary 4.2 and Theorem 4.3 from \cite{Warner}.

\subsection{Proposition} \textit{Let $k,\lambda$ as above. Define $t_0$ be the smallest positive solution $s$ of $\cot (\sqrt{k}s)=-\lambda/\sqrt{k}$. Then $$| \operatorname{Jac}_x (\exp_n)(x,t)| \leq 1 + C |t|^n + o(|t|^{n+1})$$ for $t\in[-t_0,t_0]$, where $C$ and $o(|t|^{n+1})$ depend only on $\lambda, k$ and $n$.}

\

Then by the area formula we obtain

\subsection{Corollary} \textit{Let $k,\lambda$ as above. There exists $\delta_0>0$, depending only $k$ and $\lambda$ such that $$\mathcal{H}^{n-1}(\{d_{\nu(N)} = \delta\}) \leq (1+C |\delta|^n)  \mathcal{H}^{n-1}(N)$$ for every $\delta\leq|\delta_0|$, where $C$ depends only on $\lambda, k$ and $n$.}

\subsection{Hypersurfaces parallel to level sets near a nondegenerate critical point}

In this section we are interested in estimating the area of parallel hypersurfaces to level sets of functions of the type $$f(x)=-(x_1^2 +\cdots + x_k^2)+(x_{k+1}^2+\cdots+x^{2}_{n})$$ for some $k\leq n$ and $x=(x_1,\dots,x_n)\in \R^n$. In this situation we cannot apply the results from the last section since the hypersurfaces $\{f(x)=s\}$ does not have bounded shape operator as $s\to 0$, but given that our ambient space is $\R^n$, the geometry is quite simple and we can argue by pointing out directly some geometric properties of these level sets.

First remember that whenever we have an oriented embedded hypersurface $N \subset \R^n$, with a normal Gauss map $n: N\to \R^n$, the exponential map can be written as  $\exp_n(x,t)=x+t \cdot n(x)$ and $$d_x \exp {(x,t)}= \operatorname{Id} +\ t\cdot dn(x).$$ Choosing an orthonormal basis $e_1,\cdots, e_{n-1}$ for $TN_x$ for which the matrix of the shape operator is  in diagonal form, we see that

$$| \operatorname{Jac}_x ( \exp_n){(x,t)}| = | \det((1+t\lambda_1)e_1,\dots,(1+t\lambda_{n-1})e_{n-1})| $$ $$ \leq |\Pi_{i=1}^{n-1}(1+t\lambda_i)|\leq 1 + C |t|^{n-1} | \operatorname{Jac}_x(n)(x)|,$$ where $C$ is a constant depending only on $n$.

In particular, by the same arguments as in \ref{parallelhypersurfaces}, $$\mathcal{H}^{n-1}(exp_n(N,t))\leq \mathcal{H}^{n-1}(N)+C |t|^{n-1}\int_{N} | \operatorname{Jac}_x(n)|.$$

By the area formula, $\int_{N} | \operatorname{Jac}_x(n)|$ is just the measure of the region in $\mathbb{S}^{n-1}$ covered by the Gauss map (possibly with overlaping). In the specific case of $N=f^{-1}(s)$, where $$f(x)=-(x_1^2 +\cdots + x_k^2)+(x_{k+1}^2+\cdots+x^{2}_{n}),$$ we can show that the Gauss map is injective, in particular its image has measure bounded by the volume of $\mathbb{S}^{n-1}$. 

The normal at $x$ is given by $n(x)=\nabla f(x)/|\nabla f(x)|$ where $$\nabla f(x)=2 (-x_1,\dots,-x_k,x_{k+1},\dots,x_n).$$ 

First, let $s\neq 0$ and $x,y\in N=f^{-1}(s)$ two different points, such that $n(x)=n(y)$. This implies that $\nabla f(x) \parallel \nabla f(y)$ and $x \parallel y$, i.e. there exists $\alpha\neq 0$ such that $ x = \alpha y$, but $s=f(x)=f(\alpha y) =\alpha^2 f(y)=\alpha^2 s$. Since $x\neq y$ it must be $x=-y$ and $n(x)=-n(y)$ which is a contradiction.

For $s=0$, the level set is invariant by homotheties of $\R^n$, in particular, one of its principal curvatures is always zero. Then $| \operatorname{Jac}(n)|\equiv 0$ in this case and the image of the Gauss map has measure zero.

Summarizing we obtain the following

\subsection{Lemma}\textit{Let $f$ be as in the beginning of the section. There exists a constant $C$ depending only on $n$, such that for all $s,t\in \R$ $$\mathcal{H}^{n-1}(exp_n(N\cap U,t))\leq (1+C |t|^n) \mathcal{H}^{n-1}(N\cap U),$$ where $N=f^{-1}(s)$ and $U\subset \R^n$ is a bounded open set.}

\

The following consequence of the lemma above is useful when dealing with sweepouts generated by the level sets of a Morse function. 

\subsection{Corollary}\textit{Suppose that $\psi_s \in C^{\infty}([-\e,\e]:\operatorname{Diff}(U:\R^n))$, where $U\in \R^n$ is an open set of $\R^n$. Given $\tilde{U}\subset\subset U$, define $\Sigma_s=\psi_s(f^{-1}(s)\cap \tilde{U})$. Then, there exists a constant $C$ depending only on $\psi,\tilde{U}$ and $n$ such that for all $s\in[-\e,\e]$ and $\delta\in \R$ we have $$\mathcal{H}^{n-1}(exp_{n_s}(\Sigma_s,\delta))\leq (1+C |\delta|^n) \mathcal{H}^{n-1}(\Sigma_s),$$ where $n_s$ is any Gauss map associated to $\Sigma_s$}.

\

\


\section{Appendix B: Limit-interface on manifolds}\label{manifolds.section}

We now indicate how the convergence and regularity results for phase transitions in bounded open sets of $\R^n$, can be extended to general manifolds. Since the proofs are essentially the same, the content of this section is not self-contained, but rather it is intended to serve as companion for adapting the arguments from \cite{HutchinsonTonegawa,Tonegawa,TonegawaWickramasekera}. This has been done before in the special case of closed 2-dimensional Riemannian manifold \cite{Padilla}.

\subsection{Remark}\label{remarkpde} Several arguments from \cite{HutchinsonTonegawa,TonegawaWickramasekera} involve the use of blow-up arguments and elliptic estimates, which are local in nature and can be carried out similarly in our context using normal coordinates. The Laplace-Beltrami operator coincide, in these coordinates, with the Laplacian of $\R^n$ at the origin, and usually error terms can be corrected if we restricting computations to a small neighborhood of the origin. An exception appears, for example, when trying to generalize Lemma 5.2 from \cite{HutchinsonTonegawa}, since we must deal with the error associated to the function $\psi$, but we can still control these terms with a standard application of Harnack's inequality and Schauder's estimates. 

\

Multiplying equation (\ref{eq:theequation}) by $\nabla u \cdot g$ and integrating by parts, we obtain the following useful formula. 

\begin{equation}\label{parts}
\int_{\{|\nabla u|>0\}} \big( \operatorname{div} g -  \nabla_{\nu}\  g \cdot \nu  \big) \e | \nabla u|^2 \  = \int_M \bigg( \e \frac{|\nabla u|^2}{2} - \frac{W(u)}{\e} \bigg)\operatorname{div} g,
\end{equation} 
where both integrals are with respect to the volume form of $M$; $\operatorname{div}$ denotes the divergence operator of $M$; $g$ is a smooth tangent vector field on $M$ and $\nu=\nabla u /|\nabla u|$.

\subsection{Local monotonicity formula}

One important step in \cite{HutchinsonTonegawa} is the derivation of a local monotonicity formula for the energy functional. This can be done using equation (\ref{parts}) in the following way. 

Let $\varphi:\R\to\R$ be a smooth function such that $\varphi(x)=1$ if $x\leq 0$, $\varphi(x)=1$ if $x>1$ and $0\leq \varphi(x)\leq 1$ for all $x\in \R$. Given $x\in M$, define a vector field in a neighborhood of $x$ by the formula $g=r\psi(r)\nabla r$, where $r$ is the distance to $x$ and $\psi(s)=\varphi(\frac{s-\rho}{\delta})$. Note that $\psi(r)$ converges, as $\delta\to 0$, to the characteristic function of the normal ball of radius $\rho$ centered at $x$.

In a neighborhood of $x$ we can assume that $K_M$, the sectional curvature of $M$, is bounded by $k$. An application of the Hessian comparison theorem (see Lemma 7.1 from \cite{ColdingMinicozzi}) gives
\begin{equation}\label{hessiandistance}
\bigg|(\operatorname{Hess} r )(X,X) -\frac{1}{r} \big|X - (X\cdot \nabla r )\nabla r \big|^2\bigg|\leq \sqrt{k}  
\end{equation} for $|X|=1$ and $r$ small enough.

Plugging $g$ into equation (\ref{parts}), using (\ref{hessiandistance}) and making $\delta \to 0$ we obtain
\begin{dmath}
-(n-1) \int_{B_\rho} e_\e(u) + \rho \int_{\partial B_\rho} e_\e(u) \geq \int_{B_\rho} (-\xi_\e)(u) + \e\rho\int_{\partial B_\rho}(\nabla u\ \cdot \nabla r)^2 - \sqrt{k} \int_{B_\rho}r (e_\e(u)+\e|\nabla u|^2),
\end{dmath}
where $e_\e(u)=\e|\nabla u|^2/2+W(u)/\e$ is the energy integrand and $\xi_\e(u)=\e|\nabla u|^2/2-W(u)/\e$ is the discrepancy function.

 Finally, diving by $\rho^{-n}$ and multiplying by the exponential function $e^{m\rho}$ with $m\geq 3\sqrt{k}$ we obtain the formula $$\frac{d}{d\rho} \bigg(e^{m\rho} \rho^{-n+1} \int_{B_\rho}e_\e(u) \bigg) \geq e^{m\rho} \rho^{-n+1} \int_{B_\rho}(-\xi_\e)(u),$$ the rest of the proof is the same as in \cite{HutchinsonTonegawa}.

\subsection{Stationarity and integrality}

In \cite{TonegawaPadilla} it is shown \begin{equation}\label{firstvariation} \delta V_\e(g) = \int_{\{|\nabla u|>0\}} \big( \operatorname{div} g -  \nabla_{\nu}\  g \cdot \nu  \big) \e | \nabla w|. \end{equation} 
On the other hand, in \cite{HutchinsonTonegawa}, the proof that $\xi_\e \to 0$ and $\e | \nabla u|^2 - 2|\nabla w|\to 0$ $L^1_{\text{loc}}(M)$, only involves local elliptic estimates, and the same is true in our context. This and equation (\ref{parts}) imply that the limit-interface is a stationary varifold. The rectifiability and integrality follow from the density estimates in \cite{HutchinsonTonegawa} which depend only on the local monotonicity formula and standard elliptic estimates (see Remark \ref{remarkpde}).

\subsection{Generalized second fundamental form}\label{gsff}

From now on we assume that $u_\e$ is a stable solution of equation (\ref{eq:theequation}). The following stability inequality for $u_\e$ is equivalent to the stability inequality for minimal hypersurfaces (see  \cite{Farina}, Theorem 6). \begin{equation}\label{stability} \int_M \bigg( \operatorname{Ric}(\nabla u_\e,\nabla u_\e) +|\operatorname{Hess}u_\e|^2 - \big|\nabla|\nabla u_\e| \big|^2 \bigg) \phi^2 \leq \int_M |\nabla u_\e|^2 |\nabla \phi|^2\end{equation}

For the definition of \textit{generalized second fundamental form} we refer the reader to \cite{Hutchinson,Tonegawa}. In \cite{Tonegawa} it is shown that the limit-interface have a generalized second fundamental form that satisfies a stability inequality. We adapt what is done in \cite{Hutchinson,Tonegawa}, after embedding $M$ isometrically into some $\R^p$.

For every $x\in M$,  let $P=P(x)$  be the projection onto the subspace $T_x M$ and $P_{ij}$ its coordinates on a orthonormal basis, $e_1,\dots,e_p$, of $\R^p$. Let $\nu_k$ be the coordinates of the vector $\nu=|\nabla u_\e|/|\nabla u_\e|$, whenever $|\nabla u_\e|>0$ (notice that $\nu$ depends on $\e$). We denote by $\nu\otimes\nu$ the projection onto the vector $\nu$, with coordinates $(\nu\otimes\nu)_{ij}=\nu_i\nu_j$. Then $(P-\nu\otimes\nu)$ is the projection onto the subspace orthogonal to $\nu$ in $T_x M$.

Denote by $\nabla$ and $\operatorname{div}$, the connection on $M$ (or the intrinsic gradient if applied to a function) and the divergence operator on $M$, respectively. The coordinates of $\R^n\times {Gr_{n-1}(\R^p)}$ in the basis $e_1,\dots,e_p$, we denote by $x_i$, $S_{lk}$, for $i,l,k=1,\dots,p$, and the partial derivatives in the direction of the vectors of the basis by $D_i=\partial_{x_i}$ and $D^*_{lk}$.

The definition of the second fundamental form involves functions $A_{ijk}$ (see \cite{Hutchinson}). For every $\e>0$ and $x\in M$ such that $|\nabla u^\e|>0$ define $$A^{\e}_{ijk}(x,S):=S_{si}\partial_{x_s}(P_{jk}-\nu_j\nu_k)$$ and $$(B^\e)_{ij}^k(x,S)= S_{lj} \big(A_{ikl}- S_{im} \partial_{x_m} P_{kl}\big).$$

Given $\varphi\in C^2(\R^p\times Gr_{n-1}(\R_p))$, let $g$ be the vector field tangent to $M$ defined by $$g=\phi (P-\nu\otimes\nu) (e_i),$$ where $\phi \in C^2(\R^n)$ is defined by $\phi(x)=\varphi(x,P-\nu\otimes\nu)$ (it is not true in general that $\phi \in C^2(\R^n)$, due to the singularities of $\nu$, but we can proceed as in \cite{Tonegawa}. Define $\phi_s(x)=\varphi\big(x,P-\frac{\nabla u^\e\otimes \nabla u^\e}{s+|\nabla u^\e|^2}\big)$, do all the computations with $\phi_s$ instead of $\phi$ and then make $s\to 0$).

Plugging $g$ into equation (\ref{parts}) we obtain $$\int \big( S_{si} D_s\phi + A^\e_{kik}\phi + A^\e_{ilj} D_{lj}^* \phi \big) dV^\e = \int_M \bigg( \e \frac{|\nabla u|^2}{2} - \frac{W(u)}{\e} \bigg)\operatorname{div} g$$

We assert that all the terms on this equation can be bounded as \textit{measure-function} pairs (see \cite{Hutchinson,Tonegawa}), using the stability inequality (\ref{stability}), and that, in that sense, the righthand side of this equation goes to zero as $\e\to 0$ as in \cite{Tonegawa}. Also, after passing to a subsequence, the functions $A_{ijk}^{\e}$ converge to functions $A_{ijk}$ satisfying $$\int \big( S_{si} D_s\phi + A_{kik}\phi + A_{ilj} D_{lj}^* \phi \big) dV = 0,$$ where $V$ is the limit varifold. Then $B_{ij}^k$ is the generalized fundamental form of $V$ as a $(n-1)$-varifold of $M$ (see \cite{Hutchinson,Tonegawa}). 

It is left to see that $(A_{ijk}^{\e},V^\e)$ are bounded as measure-function pairs. Notice that it is possible to compare their values in coordinates with the integrand $\big|\operatorname{Hess}u_\e|^2 - \big|\nabla|\nabla u_\e| \big|^2$ from the stability inequality. In fact, for $x\in M$, with $\nabla u_\e\neq 0$, choose $e_1,\dots,e_p$ such that $T_x M=\langle e_1,\dots,e_n \rangle=\R^n$ and $\nu=e_n$. The following observations will help make the computations easier.

Let $\overline B$ denote the second fundamental form of $M$ as a submanifold of $\R^p$ and $\overline B_{ij}^{k}=\overline B(e_i,e_j)\cdot e_k$, $\nabla u = (u_1,\cdots, u_n)$ and $\partial_{x_i}\partial_{x_j}u=u_{ij}$. Then

\begin{enumerate}
\item [{(i)}] $P_{ij}=1$ if $i=j\leq n$ and = 0 other wise.
\item [{(ii)}] $\partial_{x_s}P_{ij}=0$ if $i,j\leq n$.
\item [{(iii)}] $\partial_{x_s} \nu_n=0$ .
\item [{(iv)}] $\overline B_{ij}^k=P_{lj} P_{is} \partial_{x_s} P_{kl}.$
\item [{ (v)}] $u_{ik}=u_n \partial_{x_i} P_{kn}$ if $i\leq n$ and $k>n$. 
\end{enumerate}

{(i)} follows from the chose of the basis. {(ii)} follows differentiating the identity $P_{ij}=P_{ik}P_{kj}$ and applying {(i)}. $(iii)$ is similar.  {(iv)} is proved in \cite{Hutchinson} and {(v)} follows from {(iv)}.

Using these formulas we obtain 
$$\big|\operatorname{Hess}u_\e|^2 - \big|\nabla|\nabla u_\e| \big|^2= \sum_{i,j}^{n-1}u_{ij}^2+\sum_{j}^{n-1}u_{nj}^2.$$ and 
$$\sum_{i,j,k}^{p}|A^\e_{ijk}|^2= 2\sum_{i,j}^{n-1} \frac{{u_{ij}}^2}{u_n^2}+\sum_i^{n-1}\sum_{j,k\geq n+1}(\partial_{x_i}P_{jk})^2.$$ Notice that the last sum in the above equality depends only on $M$. Combining these with the stability formula (\ref{stability}) we obtain the boundedness of $|A^\e_{ijk}|^2$.

\subsection{Stability and regularity}

In Section 18 of \cite{Wickramasekera}, the stability of a varifold $V$ is defined by means of an inequality satisfied by the pushforward of $V$ to a tangent space of $M$ via the inverse of the exponential map. This inequality is the same found in \cite{SchoenSimon} and is a consequence of the classical stability inequality for hypersurfaces.

That the limit-interface satisfies the same inequality is a consequence of formula (\ref{stability}). In fact, proceeding as in \ref{gsff}, $$|B^\e|^2=\sum_{ijk}^{p}\bigg((B^\e)_{ij}^k \bigg)^2\leq \sum_{ij}^{n-1}\frac{u_{ij}^2}{|\nabla u|^2},$$ and we can relate this to (\ref{stability}) as in \ref{gsff}. The terms obtained in the inequality are exactly the same appearing in the stability inequality of \cite{Wickramasekera}.

The rest of the regularity proof for stable limit-interface in \cite{TonegawaWickramasekera}, involves the study of the behavior of the tangent cones of $V$. These computations are local and depend only on elliptic estimates and blowup arguments, which also hold in our context.

\bibliographystyle{siam}
\bibliography{references}

\begin{thebibliography}{10}

\bibitem{ByeonRabinowitz}
{\sc J.~Byeon and P.~H. Rabinowitz}, {\em A note on mountain pass solutions for
  a class of {A}llen-{C}ahn models}, RIMS kokyuroku, 1881 (2014), pp.~1--17.

\bibitem{CaffarelliVasseur}
{\sc L.~Caffarelli and A.~Vasseur}, {\em The {D}e {G}iorgi method for
  regularity of solutions of elliptic equations and its applications to fluid
  dynamics}, Discrete Contin. Dyn. Syst. Ser. S, 3 (2010), pp.~409--427.

\bibitem{AllenCahn}
{\sc J.~Cahn and S.~Allen}, {\em A microscopic theory for domain wall motion
  and its experimental verification in {Fe}-{Al} alloy domain growth kinetics},
  Le Journal de Physique Colloques, 38 (1977), pp.~C7--51.

\bibitem{ColdingDelellis}
{\sc T.~H. Colding and C.~De~Lellis}, {\em The min-max construction of minimal
  surfaces}, in Surveys in differential geometry, {V}ol.\ {VIII} ({B}oston,
  {MA}, 2002), Surv. Differ. Geom., VIII, Int. Press, Somerville, MA, 2003,
  pp.~75--107.

\bibitem{ColdingMinicozzi}
{\sc T.~H. Colding and W.~P. Minicozzi}, {\em A course in minimal surfaces},
  vol.~121, American Mathematical Soc., 2011.

\bibitem{RosenbergHyp}
{\sc P.~Collin, L.~Hauswirth, L.~Mazet, and H.~Rosenberg}, {\em Minimal
  surfaces in finite volume non compact hyperbolic $3 $-manifolds}, arXiv
  preprint arXiv:1405.1324,  (2014).

\bibitem{DelellisPellandini}
{\sc C.~De~Lellis and F.~Pellandini}, {\em Genus bounds for minimal surfaces
  arising from min-max constructions}, Journal f{\"u}r die reine und angewandte
  Mathematik (Crelles Journal), 2010 (2010), pp.~47--99.

\bibitem{DelellisTasnady}
{\sc C.~De~Lellis and D.~Tasnady}, {\em The existence of embedded minimal
  hypersurfaces}, Journal of Differential Geometry, 95 (2013), pp.~355--388.

\bibitem{Padilla}
{\sc H.~del Rio~Guerra, C.~Garza-Hume, and P.~Padilla}, {\em Geodesics, soap
  bubbles and pattern formation in riemannian surfaces}, The Journal of
  Geometric Analysis, 13 (2003), pp.~595--604.

\bibitem{EvansGariepy}
{\sc L.~C. Evans and R.~F. Gariepy}, {\em Measure theory and fine properties of
  functions}, vol.~5, CRC press, 1991.

\bibitem{Farina}
{\sc A.~Farina, Y.~Sire, and E.~Valdinoci}, {\em Stable solutions of elliptic
  equations on {R}iemannian manifolds}, Journal of Geometric Analysis, 23
  (2013), pp.~1158--1172.

\bibitem{Ghoussoub}
{\sc N.~Ghoussoub}, {\em Duality and perturbation methods in critical point
  theory}, vol.~107, Cambridge University Press, 1993.

\bibitem{HeintzeKarcher}
{\sc E.~Heintze and H.~Karcher}, {\em A general comparison theorem with
  applications to volume estimates for submanifolds}, Ann. Sci. Ecole Norm.
  Sup, 11 (1978), pp.~451--470.

\bibitem{Hutchinson}
{\sc J.~Hutchinson}, {\em Second fundamental form for varifolds and the
  existence of surfaces minimising curvature}, Indiana Univ. Math. J., 35
  (1986), pp.~45--71.

\bibitem{HutchinsonTonegawa}
{\sc J.~E. Hutchinson and Y.~Tonegawa}, {\em Convergence of phase interfaces in
  the van der {W}aals-{C}ahn-{H}illiard theory}, Calculus of Variations and
  Partial Differential Equations, 10 (2000), pp.~49--84.

\bibitem{Ilmanen}
{\sc T.~Ilmanen}, {\em Convergence of the {A}llen-{C}ahn equation to {B}rakke's
  motion by mean curvature}, J. Differential Geom, 38 (1993), pp.~417--461.

\bibitem{Ketover}
{\sc D.~Ketover}, {\em Degeneration of min-max sequences in 3-manifolds}, arXiv
  preprint arXiv:1312.2666,  (2013).

\bibitem{Li}
{\sc M.~M.-c. Li}, {\em A general existence theorem for embedded minimal
  surfaces with free boundary}, Communications on Pure and Applied Mathematics,
  68 (2015), pp.~286--331.

\bibitem{MarquesNevesDuke}
{\sc F.~C. Marques and A.~Neves}, {\em Rigidity of min-max minimal spheres in
  three-manifolds}, Duke Mathematical Journal, 161 (2012), pp.~2725--2752.

\bibitem{MarquesNevesInfinite}
{\sc F.~C. Marques and A.~Neves}, {\em Existence of infinitely many minimal
  hypersurfaces in positive {R}icci curvature}, arXiv preprint arXiv:1311.6501,
   (2013).

\bibitem{MarquesNevesWillmore}
{\sc F.~C. Marques and A.~Neves}, {\em Min-max theory and the {W}illmore
  conjecture}, Ann. of Math. (2), 179 (2014), pp.~683--782.

\bibitem{Rosenberg}
{\sc L.~Mazet and H.~Rosenberg}, {\em Minimal hypersurfaces of least area},
  arXiv preprint arXiv:1503.02938,  (2015).

\bibitem{MizunoTonegawa}
{\sc M.~Mizuno and Y.~Tonegawa}, {\em Convergence of the {A}llen--{C}ahn
  equation with {N}eumann boundary conditions}, arXiv preprint arXiv:1403.5624,
   (2014).

\bibitem{Modica}
{\sc L.~Modica}, {\em The gradient theory of phase transitions and the minimal
  interface criterion}, Archive for Rational Mechanics and Analysis, 98 (1987),
  pp.~123--142.

\bibitem{Montezuma}
{\sc R.~Montezuma}, {\em Min-max minimal hypersurfaces in non-compact
  manifolds}, arXiv preprint arXiv:1405.3712,  (2014).

\bibitem{Pacard}
{\sc F.~Pacard}, {\em The role of minimal surfaces in the study of the
  {A}llen-{C}ahn equation}, Geometric Analysis: Partial Differential Equations
  and Surfaces: UIMP-RSME Santal{\'o} Summer School Geometric Analysis, June
  28-July 2, 2010, University of Granada, Granada, Spain, 570 (2012), p.~137.

\bibitem{TonegawaPadilla}
{\sc P.~Padilla and Y.~Tonegawa}, {\em On the convergence of stable phase
  transitions}, Communications on pure and applied mathematics, 51 (1998),
  pp.~551--579.

\bibitem{PisanteAllen13}
{\sc A.~Pisante and F.~Punzo}, {\em {A}llen-{C}ahn approximation of mean
  curvature flow in {R}iemannian manifolds i, uniform estimates}, arXiv
  preprint arXiv:1308.0569,  (2013).

\bibitem{Pitts}
{\sc J.~T. Pitts}, {\em Existence and regularity of minimal surfaces on
  Riemannian manifolds}, no.~27 in Mathematical Notes, Princeton University
  Press, Princeton, 1981.

\bibitem{RabinowitzBook}
{\sc P.~H. Rabinowitz}, {\em Minimax methods in critical point theory with
  applications to differential equations}, no.~65 in Conference Board of the
  Mathematical Science, American Mathematical Soc., 1986.

\bibitem{Savin}
{\sc O.~Savin}, {\em Phase transitions, minimal surfaces and a conjecture of
  {D}e {G}iorgi}, Current developments in mathematics,  (2009), pp.~59--113.

\bibitem{SchoenSimon}
{\sc R.~Schoen and L.~Simon}, {\em Regularity of stable minimal hypersurfaces},
  Communications on Pure and Applied Mathematics, 34 (1981), pp.~741--797.

\bibitem{SchoenSimonYau}
{\sc R.~Schoen, L.~Simon, and S.-T. Yau}, {\em Curvature estimates for minimal
  hypersurfaces}, Acta Mathematica, 134 (1975), pp.~275--288.

\bibitem{SimonBook}
{\sc L.~Simon}, {\em Lectures on geometric measure theory}, Centre for
  Mathematical Analysis, Australian National University Canberra, 1984.

\bibitem{SimonSmith}
{\sc L.~Simon and F.~Smith}, {\em On the existence of embedded minimal
  2-spheres in the 3-sphere, endowed with an arbitrary metric}, preprint.

\bibitem{Sternberg}
{\sc P.~Sternberg}, {\em The effect of a singular perturbation on nonconvex
  variational problems}, Archive for Rational Mechanics and Analysis, 101
  (1988), pp.~209--260.

\bibitem{Tonegawa}
{\sc Y.~Tonegawa}, {\em On stable critical points for a singular perturbation
  problem}, Communications in Analysis and Geometry, 13 (2005), pp.~439--459.

\bibitem{TonegawaSurvey}
\leavevmode\vrule height 2pt depth -1.6pt width 23pt, {\em Applications of
  geometric measure theory to two-phase separation problems}, Sugaku
  Expositions, 21 (2008), p.~97.

\bibitem{TonegawaWickramasekera}
{\sc Y.~Tonegawa and N.~Wickramasekera}, {\em Stable phase interfaces in the
  van der {W}aals--{C}ahn--{H}illiard theory}, Journal f{\"u}r die reine und
  angewandte Mathematik (Crelles Journal), 2012 (2012), pp.~191--210.

\bibitem{Warner}
{\sc F.~Warner}, {\em Extension of the {R}auch comparison theorem to
  submanifolds}, Transactions of the American Mathematical Society,  (1966),
  pp.~341--356.

\bibitem{Wickramasekera}
{\sc N.~Wickramasekera}, {\em A general regularity theory for stable
  codimension 1 integral varifolds}, Annals of Mathematics, 179 (2014),
  pp.~843--1007.

\bibitem{Zhou}
{\sc X.~Zhou}, {\em Min-max minimal hypersurface in $({M}^{n+1}, g)$ with
  $\operatorname{Ric}_g > 0$ and $2\leq n\leq 6$}, arXiv preprint
  arXiv:1210.2112,  (2012), pp.~22460--320.

\bibitem{Zhou2015}
{\sc X.~Zhou}, {\em Min-max hypersurface in manifold of positive ricci
  curvature}, arXiv preprint arXiv:1504.00966,  (2015).

\end{thebibliography}

\end{document}